%% file: ms.tex
\documentclass{amsart}
\usepackage{amssymb}
\usepackage{amsmath}
\usepackage[hyperref]{degt}
\def\itemrefform#1{$(#1)$}

\addtheorem{Algorithm}{algorithm}

\def\2{\color{red}}
\def\3{\color{blue}}
\def\4{\color{magenta}}
\def\Dg:{{\bf Dg:\enspace}\ignorespaces}

\def\KK-{\pdfstr{K3}{$K3$}-\penalty0\ignorespaces}

\let\0\varnothing

\let\graph\Gamma
\let\pencil\Pi
\let\fiber\Sigma
\let\Weyl\Delta
\def\fp{\Weyl^{\sharp}}

\def\Fano{\Cal{F}}
\def\Cone{\Cal{C}}
\def\Fn{\operatorname{Fn}}
\def\extended{^{\mbox{\tiny ex}}}
\def\extended{^\text{\rm ex}}
\def\Fex{\Fn\extended}
\def\rt{\operatorname{\frak{rt}}}
\def\root{\operatorname{root}}
\def\sat{\operatorname{sat}}
\def\satex{\sat\extended}
\def\girth{\operatorname{girth}}
\def\val{\operatorname{val}}

\def\spp{\operatorname{sp}}
\def\NS{\operatorname{NS}}

\def\hp{h}

\def\graphl{\operatorname{Fn}}
\def\graphe{\graphl\extended}

\def\qgen{spanned }

\def\resxf{\smash{\tilde{X}_4}}
\def\hhp{\pi^*{\mathcal O}_{X_4}(1)}

\def\thmcit{Theorem}

\def\AAmple{Abundant}
\def\aample{abundant}

\def\resp.{resp\PERIOD}
\let\geq\ge

\let\M=M
\let\o=o

\let\ex=e
\let\iso=p

\def\basis#1{\{#1\}}

\def\CK{\Cal K}

\def\GK{\frak G}
\def\bL{\bold{L}}

\def\CC{\Cal C}

\def\CS{\Cal S}
\def\CG{\Cal G}

\makeatletter
\def\newconfig#1#2{\expandafter\gdef\csname config=#1\endcsname{#2}}

\def\theconfig#1{\hyperlink{config=#1}{\csname config=#1\endcsname}}
\def\confighook#1#2{\smash{\raise9pt\hbox{\hypertarget{config=#1}{}}}%
 \protected@write\@auxout{}{\string\newconfig{#1}{#2}}%
 \csname config=#1\endcsname}
\pdef\aconfig#1#2{\ifmmode\confighook{#1}{#2}\else$\confighook{#1}{#2}$\fi}
\pdef\config#1{\ifmmode\theconfig{#1}\else$\theconfig{#1}$\fi}
\makeatother

\makeatletter
\let\HyPsd@CatcodeWarning\@gobble
\makeatother

\def\ta{\tilde{a}}

\let\trig\vartriangle
\let\set\Theta

\def\val{\operatorname{val}}
\def\girth{\operatorname{girth}}
\def\Pencils{\Cal{P}}
\def\Sets{\Cal{T}}
\def\Sec{\bold{A}}
\let\cp=\pi
\let\cs=\theta
\def\smooth{^\mathrm{s}}

\let\less=\prec
\let\lesseq=\preccurlyeq

\let\grteq=\succcurlyeq
\let\trig\vartriangle
\let\icup\sqcup

\def\pert{\mathrel\vartriangleleft}

\def\comp{\mathrel\leftthreetimes}
\let\comp\vdash

\def\bu{\bold u}
\def\bv{\bold v}
\def\gm{\frak m}
\def\CS{\Cal{S}}
\def\gS{\frak S}
\def\fixed{\underline{\text{fixed}}}
\def\bQ{\bold Q}
\def\bX{\bold X}
\def\bY{\bold Y}
\def\bZ{\bold Z}
\def\Sym{\operatorname{Sym}}

\title{Lines on \KK-quartics via triangular sets}

\author{Alex Degtyarev}

\address{%
	Department of Mathematics\\
	Bilkent University\\
	06800 Ankara, TURKEY}

\email{
	degt@fen.bilkent.edu.tr}

\author{S\l awomir Rams}
\address{Institute of Mathematics, Jagiellonian University,
	ul. {\L}ojasiewicza 6,  30-348 Krak\'{o}w, Poland}
\email{slawomir.rams@uj.edu.pl}

\thanks{%
	A.D. was partially supported by the T\"{U}B\DOTaccent{I}TAK grant 118F413.
    S.R. was  partially supported by the National Science Centre, Poland, Opus  grant
	no. 2018/31/B/ST1/02857} 



\begin{document}

	\begin{abstract}
	 We prove the  sharp upper bound of at most 52 lines on a
complex \KK-surface of degree $4$ with a non-empty singular locus.
We also classify the configurations of more than $48$ lines on smooth  complex quartics.
	\end{abstract}
	
	\maketitle

\section{Introduction}\label{S.intro}


Our main goal
is to
present an approach to study large line configurations
on complex projective \KK-quartics.
In particular,
we prove the following theorem.

\theorem[see \autoref{proof.main}]\label{thm-main}
Let $X_{4} \subset \Cp{3}(\C)$ be a degree $4$ \KK-surface  with
non-empty singular locus.
Then $X_4$ contains at most $52$ lines. Moreover,
each  \KK-quartic with at least $49$ lines contains four coplanar lines.
\endtheorem

The above  bound is sharp: the existence of a complex \KK-quartic with $52$
lines and non-empty singular locus
(two simple nodes) was shown by the first named author in
2016 (\via\ Torelli's theorem, see \cite[Theorem~1.10]{degt:singular.K3}) and
the equation of the surface in question was found by D.~Veniani, see
\cite[Example~5.5]{Veniani:char3}.

We conjecture that the quartic surface discovered in
\cite[Theorem~1.10]{degt:singular.K3} is the
only
quartic that attains the
bound of \autoref{thm-main}, but the proof of this fact is beyond the scope
of this paper.

It is well-known that the complexity of large line configurations on
projective \KK-surfaces
decreases
as the degree~$d$ of the polarization grows. In particular, a complete
classification of close to maximal configurations is known for octics (see
\cite[Theorem 1.1]{Degtyarev.Rams}) and sextics (to appear
in~\cite{Degtyarev.Rams.sextics}):
the respective upper bounds are $32$ and $36$ in the presence of a
singularity \vs. $36$ and $42$ in the smooth case.
In contrast, even though quartic surfaces
with singular points have been a subject of intensive study
ever
since the
19-th century (see, \eg, the classical treatise \cite{Jessop:quartic}),
hardly anything is known about large line configurations on such surfaces.
The main
reason
is the existence of the so-called
\emph{triangular} configurations on
quartics (see \autoref{s.classification} for the definition) --- a property
that drastically increases the complexity of the problem.  Here, we circumvent
this difficulty with the help of the so-called \emph{triangular sets}
introduced in \autoref{S.trig}.

One can easily check that the degree-$d$  Fermat surface (over $\C$) contains exactly $3d^2$  lines for $d>2$.
Moreover,  for almost all integers $d$ the Fermat surface is the best known example of a smooth complex projective surface with many lines, and the question whether smooth degree-$d$ surfaces with more lines exist remains open.
To
illustrate the power of our approach, we refine the results of \cite{DIS}
and classify all configurations of
at least $49$ lines on \emph{smooth}
quartics (\ie, the configurations that are larger than the one on the Fermat quartic).
 Remarkably, compared to \cite{DIS}, we found but three
new configurations: one of
 rank $20$
 ($\bQ_{52}'''$ previously found
in~\cite{degt:singular.K3}) and two of
rank $19$
(designated as $(52)$
and $(50)$ in \autoref{tab.smooth}). On the other hand, there are at least
$28$ configurations of $48$ lines on smooth quartics,
giving yet another reason
why $48$
 is a
reasonable threshold (\cf. also \autoref{prop.quad}
and \autoref{rem.c3} below).

\theorem[see \autoref{proof.smooth}]\label{th.smooth}
Up to isomorphism, there are $26$ configurations of at least $49$ lines on
smooth quartic surfaces, see \autoref{tab.smooth}.
They are realized by $34$ singular \rom(\latin{aka} projectively rigid\rom) surfaces
\rom($18$ real and $8$ pairs of complex conjugate\rom) and five
connected $1$-parameter families.
\table
\caption{Smooth complex quartics with at least $49$ lines}\label{tab.smooth}
\hbox to\hsize{\hss\vbox{\offinterlineskip%
\def\config#1(#2){#1}%
\def\ps#1{}%
\def\aut#1{#1}%
\def\symplectic(#1,#2){\def\symsize{#1}\def\symindex{#2}}%
\def\group#1(#2){\symsize\rlap{$_{\symindex}^{\ifnum#2>1 #2\fi}$}}%
\def\group#1(#2){(\symsize,\symindex)\ifnum#2>1 \rlap{$^{#2}$}\fi}%
\def\mat#1{#1}%
\def\mmat#1{#1}%
\def\*{\llap{$^*$}}%
\def\bU(#1){\bold{U}_{#1}}%
\def\bU(#1){\bold{U}(#1)}%
\let\+\oplus
\def\rc#1{#1}%
\def\i{'}\def\ii{''}\def\iii{'''}\def\0{}\def\v{^{\mathrm{v}}}%
\def\(#1){(#1)}%
\halign{\strut\quad\hss$#$\hss\quad&#&#&#&\hss$#$\hss\quad&\hss$#$\hss\quad&
 \ \hss$#$\hss\quad&\ $#$\hss\quad\cr
\noalign{\hrule\vskip2pt}
\graph&&&&\ls|\Aut\graph|&\ls|\Sym X_4|&(r,c)&\NS(X_4)^\perp\cr
\noalign{\vskip1pt\hrule\vskip2pt}
\input tab_smooth

\noalign{\vskip1pt\hrule}%
\crcr}}\hss}%
\endtable
\endtheorem

As a consequence,
we answer a question
left open in \cite[Addendum 1.4]{DIS}.

\addendum[see \autoref{proof.real}]\label{add.real}
The complete list of values taken by the number of \emph{real} lines on a
\emph{real} smooth quartic is
\[*
\{0,1,\ldots,49,50,52,56\}.
\]
The configurations of more than $48$ real lines on a real
smooth quartic are those
marked with a $^*$ in \autoref{tab.smooth}.
\endaddendum

Listed in \autoref{tab.smooth} are:
\roster*
\item
the name of the configuration~$\graph$ (mostly following~\cite{DIS});
the subscript always refers to the number of lines (vertices of~$\graph$);
\item
the size of the group $\Aut\graph$ of abstract graph automorphisms of~$\graph$;
\item
the group $\Sym X_4$ of simplectic automorphisms of a generic
quartic~$X_4$ with
the given Fano graph,
in the form $(\text{size},\text{index})$, referring to the
\texttt{SmallGroup} library in \GAP~\cite{GAP4};
the superscript is the index of $\Sym X_4$ in the full group $\Aut(X_4,h)$ of
projective automorphisms of~$X_4$ (if greater than~$1$);
\item
the numbers $(r,c)$ of, respectively, real and pairs of complex conjugate
components of the equilinear moduli space;
\item
the (generic, if $\rank\ge3$) transcendental lattice $T:=\NS(X_4)^\perp$;
it is marked with a $^*$ if the corresponding deformation family has
a real quartic with all lines real (see \cite[Lemma 3.8]{DIS}).
\endroster
If $T$ is not determined by~$\graph$,
each lattice is listed in a separate row (following the main entry),
and the numbers $(r,c)$ of components are itemized accordingly.

As in \cite{degt:lines},
we use the following notation for common integral lattices:
\roster*
\item
$[a]:=\Z u$ is the  lattice of rank~$1$ given by the condition $u^2=a$;
\item
$[a,b,c]:=\Z u+\Z v$, $u^2=a$, $u\cdot v=b$, $v^2=c$, is a lattice of
rank~$2$; when it is positive definite, we assume that $0<a\le c$ and
$0\le2b\le a$: then, $u$ is a shortest vector, $v$ is a next shortest one,
and the triple $(a,b,c)$ is unique;
\item
$\bU:=[0,1,0]$ is 
the unimodular even lattice of rank~$2$;
\item $L(n)$ denotes the lattice obtained by the scaling of
a given lattice~$L$ by a fixed integer  $n\in\Z$.
\endroster
In general,
we maintain the standard notation for various objects
associated to a lattice $L$ (the determinant, discriminant group, \etc.)
---see, \eg,~\cite{Conway.Sloane,Nikulin:forms}.
The
inertia indices of the quadratic form
$L \otimes \R$ are denoted by $\sigma_{\pm,0}(L)$.


\subsection{Contents of the paper}\label{s.contents}

Roughly,
the paper consists of two parts: the discrete
one
(\autoref{S.techprel}--\autoref{S.other})
and the geometric one  (\autoref{S.Proof}).

Our approach is a refinement of the technique
developed
in \cite{DIS,Degtyarev.Rams},
and
we recall the necessary facts and introduce certain technical terms
(\eg, acceptable graphs)
 in  \autoref{S.techprel}. Then,
 in \autoref{S.trig},
  we define the main technical tool,
\viz.
the triangular set, and discuss methods  of extending
a given
graph by a
collection of triangular sets. Finally,
after those preparations, we present the proof of
 \autoref{prop.trig},
which is the discrete counterpart of
the most difficult case of \autoref{thm-main}.

\autoref{S.smooth}
is a digression:
we restrict our attention to the case of smooth lattices
(\ie, we assume that the lattice contains no exceptional divisors) and
apply triangular sets to
classify geometric Fano graphs with at least  $49$ vertices.


In \autoref{S.quad} and \autoref{S.other},
we turn back to the general case (with exceptional divisors allowed)
and study the properties of triangular free Fano graphs.

Finally, in \autoref{S.Proof} we recall the definition of the Fano graph
(resp.\ extended Fano graph) of a surface
and
its relation to the geometricity of the
Fano graph of a lattice, see \autoref{th.K3quartics}
(resp.\
\autoref{th.K3quartics.ext}) and prove the
principal
results of the paper,
\viz.
\autoref{thm-main} and \autoref{th.smooth}.


\subsection{History of the problem}\label{s.sart}
As mentioned,
configurations of lines (or, more generally, smooth rational curves) on
quartic surfaces in
$\Cp{3}$ have been a subject of intensive study
ever since the 19-th century. Still, the methods of Italian school were not
efficient enough to deal with the classification of large line configurations
on such surfaces. It was
not until
the last decade that the theory of elliptic
fibrations, Mordell--Weil groups,   Torelli's theorem and progress in
algorithmic methods in the theory of lattices led to a substantial progress
in the case of smooth quartics: sharp bound for the number of lines over
fields of characteristic $p \neq 2,3$ (see \cite{Segre,rams.schuett,DIS}),
$p=3$ (see \cite{rams.schuett:char3}), $p=2$ (see \cite{degt:supersingular}),
the classification of large configurations (see \cite{DIS}), explicit
equations of surfaces with many lines (see \cite{Veniani:equations.publ} and
the bibliography therein). Strangely enough, the  Bogomolov--Miyaoka--Yau
inequality yields no bounds 
 in the case of quartics (see \cite{Miyaoka:lines}). 	

In contrast, in spite
of long interest (see, \eg, the classical text \cite{Jessop:quartic}), far
less is known in the case of quartic surfaces with singular points ---
essentially,
it was
only shown that,
over a field of characteristic $p \neq 2$,
the
number of lines on a quartic with singularities cannot exceed the maximal
number of lines on a smooth quartics
--- see
\cite{Veniani,Veniani:char3,Rams.Gonzales}.  For $p=2$ the
maximal number of lines on a quartic with singularities is $68$
(\vs. $60$ in the smooth case, see
\cite{Veniani:char2}) and we do know projective models of surfaces that attain
this maximum (see \cite{rams.schuett:char3}).

A refinement of the method pioneered in \cite{DIS} led to the complete
picture of large line configurations on smooth degree-$d$ \KK-surfaces for
$d>2$ in \cite{degt:lines}. Vinberg's algorithms combined with the above
methods yield a means to  classify the large configurations of lines on
degree-$d$ \KK-surfaces with at worst Du Val singularities for $d>4$ (see
\cite{Degtyarev.Rams,Degtyarev.Rams.sextics}). The methods of
\cite{Degtyarev.Rams} are not sufficient to deal with the case of quartics
(\ie, $d=4$):  
the existence of triangles (\ie, \smash{$\tA_{2}$}-configurations) of lines
and the fact that,
on quartic surfaces,
said triangles may interlace lead to
numerous configurations that are excluded on degree-$d$ \KK-surfaces for
$d>4$.
In the present paper, we discuss  
an approach to deal with such
configurations.
However, in order to keep our exposition compact,  we apply our
method to find the maximal number of lines on a complex \KK-quartics with
non-empty singular locus,
but we do not try to classify all
configurations of $52$ lines.


\subsection{Acknowledgements}

This paper was mostly written during our research stay at the
\emph{Max-Planck-Institut f\"{u}r Mathematik}, Bonn. We are grateful to MPIM for
creating perfect working conditions. S.R. thanks IM PAN
(Cracow, Poland) for the support that enabled him to complete this project.

\section{Preliminaries}\label{S.techprel}

In this section we recall the main technical tools
that we use in our work.
To shorten the exposition,
we focus on the case of $4$-polarized $2$-admissible lattices and graphs.
The details and more general statements can be found in \cite{Degtyarev.Rams}.
To keep the exposition continuous,
we assume the reader familiar with the
basics of the theory of $K3$-surfaces, $(-2)$-curves, \etc. and adopt a
formal graph-theoretical language. The relation of graphs
considered in \autoref{S.techprel}--\autoref{S.other} to
the problem at hand, \ie, lines on
quartic surfaces, is briefly discussed in \autoref{s.linesonquartics} below,
right before the proofs of the principal results of the paper.

\subsection{Polarized lattices} \label{s.lattices}
Recall that a nondegenerate lattice~$S$ is called \emph{hyperbolic} if
$\Gs_+S=1$. A \emph{polarized lattice} $S\ni h$ is a hyperbolic lattice~$S$
equipped with a distinguished vector $h$ of positive square; the square $h^2$
is called the \emph{degree} of the polarization and $S$ is said to be
$h^2$-polarized. Here we assume $h^2 = 4$, so whenever we speak of a
\emph{polarized lattice} we mean a \emph{$4$-polarized lattice}.
Furthermore, we confine ourselves to lines and exceptional divisors
(resp.\ only lines in \autoref{S.smooth}),
leaving out smooth rational curves of higher degree.

\remark\label{rem.index}
We make frequent use of the following obvious observation: if $S$ is
a hyperbolic lattice,
then any sublattice $N\subset S$ is either semidefinite (and then one has
$\rank\ker N=1$) or nondegenerate.
\endremark

For a polarized lattice and $n = 0,1$, one defines the sets
\[*
\root_n(S,h):=\bigl\{r\in S\bigm|\text{$r^2=-2$, $r\cdot h=n$}\bigr\}  .
\]
As in \cite[\S\,2.2]{Degtyarev.Rams},
we put $\rt(S,h)\subset h^\perp\subset S$ (resp.\ $\Cone^+(S,h)$)
to denote the sublattice generated by
$\root_0(S,h)$ (resp.\ the positive cone).
Recall that every connected component  $\fp$ of
\[*
\Cone^+(S,h) \sminus \bigcup\limits_{r^2 = -2} r^{\perp}
\]
is  a \emph{fundamental polyhedron} for the group generated by
reflections of~$S$.

By definition,  $\rt(S,h)$ is a root lattice and each
fixed Weyl chamber~$\Weyl$
for (the group generated by reflections of) $\rt(S,h)$ gives rise
to a distinguished
fundamental polyhedron $\fp$.
 We put  $\basis\Weyl$ to denote the
"outward" roots orthogonal to the walls of $\Weyl$ and define
 the (plain) \emph{Fano graph} of the polarized lattice $(S,h)$ with a distinguished
Weyl chamber $\Weyl$ for $\rt(S,h)$  as the set of vertices
\[ \label{eq-def-Fano}
\Fn_\Weyl(S,h):=\basis\fp_1:=\big\{l\in\root_1(S,h)\bigm|
\text{$l\cdot\ex\ge0$ for all $\ex\in\basis\Weyl$}\bigr\},
\]
with two vertices $l_1\ne l_2$ connected by an edge of multiplicity
$l_1\cdot l_2$. The bi-colored \emph{extended Fano graph} is defined as
\[ \label{eq-def-exFano-lattice}
\Fex_\Weyl(S,h):=\basis\fp_1\cup\basis\Weyl,
\]
with the same convention about the multiplicities of the edges and vertices~$v$ colored
according to the value $v\cdot h\in\{0,1\}$.

\definition \label{def.separating.roots}
Let  $S\ni h$ be a polarized lattice and let  $\graph$ be a subset of $\root_1(S,h)$.
\roster
\item
A Weyl chamber~$\Weyl$ is called \emph{compatible} with
$\graph$
 if $\graph\subset\Fn_\Weyl(S,h)$.

\item
A root $r\in\root_0(S,h)$ is called \emph{separating}
with respect to
$\graph$  if there is
a pair of
vertices $u,v\in\graph$ \emph{separated} by~$r$,
so that
$r\cdot u>0$ and $r\cdot v<0$.
\endroster
\enddefinition

Finally, in order to use general theory of \KK-surfaces
 in the sequel we need the following definition.

\definition \label{def.lattice.admissible}
A polarized lattice $S\ni h$ is called:
\roster
\item
\emph{admissible}, if
there is no vector  $\iso \in S$  such that
$\iso^2=0$ and
$\iso\cdot h=2$;
\item
\emph{geometric},
if it is admissible and there exists a primitive
isometry
\[*
S\into\bL:=2\bE_8\oplus3\bU.
\]
\endroster
\enddefinition

\subsection{Subgeometric and geometric graphs}  \label{s.plain.graphs}

Let $\graph$ be a (plain) graph. To~$\graph$
 we associate
 the polarized lattice
\[ \label{eq:fanolattice}
\Fano(\graph):=(\Z\graph+\Z h)/\ker,\qquad h^2=4,\quad
\text{$h\cdot v=1$ for $v\in\graph$}.
\]
where $\Z\graph$ is the lattice freely generated by the
vertices $v\in\graph$,
so that
$u\cdot v=n$ when  $u \neq v$  are
connected by an $n$-fold edge, and $v^2=-2$ for each  $v\in\graph$.

\convention\label{conv.graphs}
 As in 
\cite[\S\,4]{Degtyarev.Rams}, we speak of \emph{polarized} graphs $\graph$
(omitting the degree which is fixed to equal~$4$),
and we apply to $\graph$
 the lattice theoretic terminology
such as
 the \emph{rank}
 $\rank\graph:=\rank\Fano{}(\graph)$ \etc.
Furthermore, we treat the vertices of~$\graph$ as vectors in $\Fano(\graph)$:
\eg, $u\cdot v\in\Z$ stands for the multiplicity of the edge $[u,v]$,
and we say that $u,v\in\graph$ \emph{intersect} if $u\cdot v=1$. The only
exception from this rule is the classification of graphs according to the
inertia indices of $\Z\graph$ rather than $\Fano(\graph)$ (which is
\emph{always} assumed hyperbolic): thus, we say that $\graph$ is
\roster*
\item
\emph{elliptic}, if $\Gs_+(\Z\graph)=\Gs_0(\Z\graph)=0$,
\item
\emph{parabolic}, if $\Gs_+(\Z\graph)=0$ and $\Gs_0(\Z\graph)>0$, and
\item
\emph{hyperbolic}, if $\Gs_+(\Z\graph)=1$ (no assumption on $\Gs_0$).
\endroster
Recall that any connected elliptic (resp.\ parabolic) graph is a Dynkin
diagram (resp.\ affine Dynkin diagram);
as in~\cite{degt:lines}, we order the isomorphism classes of affine
Dynkin diagrams according to their Milnor number, followed by $\bA<\bD<\bE$.
Recall also that, for each connected parabolic subgraph~$\fiber$, there is a
unique positive minimal generator $\kappa_\fiber\in\ker\Z\fiber$;
it has the form $\kappa_\fiber=\sum m_cc$, $c\in\fiber$, with all
$m_c>0$.



We define the \emph{perturbation order}
on the set of (isomorphism classes of) elliptic and parabolic graphs:
$\graph'\pert\graph''$ if $\graph'$ is isomorphic to an
induced subgraph of $\graph''$.
\endconvention

 Given  an isotropic subgroup
 $\CK\subset\discr\Fano(\graph)$ (\latin{aka} kernel), we consider the
finite index extension
 $\Fano(\graph,\CK)$
  of $\Fano(\graph)$ by~$\CK$
(\cf.~\cite{Nikulin:forms}).
The
pair $(\graph,\CK)$  is
  said to be \emph{extensible} if it admits a compatible Weyl chamber~$\Weyl$ for
  $\rt(\Fano(\graph,\CK),h)$ (see \autoref{def.separating.roots}):
\[*
 \graph\subset\Fn_\Weyl\Fano(\graph,\CK).
\]
Recall that, by \cite[Lemma~3.4]{Degtyarev.Rams}, we have
\begin{equation} \label{eq.l34}
\text{$(\graph,\CK)$ is extensible if and only if
 $\Fano(\graph,\CK)$ has no separating roots}
\end{equation}
(see \autoref{def.separating.roots});
moreover,
if this is the case, the compatible Weyl chamber 
  $\Weyl\subset\rt\Fano(\graph,\CK)$ is unique.
 Therefore, for an extensible pair $(\graph,\CK)$ one can define its \emph{saturation}
 and \emph{extended saturation}
 \[*
 \sat(\graph,\CK):=\Fn_\Weyl\Fano(\graph,\CK),\qquad
 \satex_{}(\graph,\CK):=\Fex_\Weyl\Fano(\graph,\CK).
 \]

A graph $\graph$ (resp.\ pair $(\graph,\CK)$)
 is called  \emph{admissible}
  if it is extensible and the lattice $\Fano(\graph)$ (\resp.\ $\Fano(\graph,\CK)$) is admissible.
 Then, an isotropic subgroup
 $\CK\subset\discr\Fano(\graph)$
 is called a \emph{geometric kernel} if the lattice
 $\Fano(\graph,\CK)$ is geometric.
 We  follow  \cite{Degtyarev.Rams} and put
 \[*
 \GK(\graph)  := \bigl\{ \CK\subset\discr\Fano(\graph)\bigm|
 \text{$\CK$ is geometric} \bigr\}.
 \]
  After those preparations we can recall the following definition. 

 \definition\label{def.madmissible}
 Let $\graph$ be a graph.
 \roster
 \item
 We call  $\graph$  \emph{subgeometric} if
 the set  $\GK(\graph)$ is non-empty.
 \item\label{graph.geometric}
 A subgeometric graph $\graph$ is called  \emph{geometric} if
 $\graph\cong\sat(\graph,\CK)$
 for a
 certain kernel $\CK \in \GK(\graph')$.
 \item
 A \emph{bi-colored} graph~$\graph'$ is
 \emph{geometric} 
  if $ \graph'\cong\Fex_\Weyl(S,h)$
for some geometric polarized lattice $S\ni h$.
\endroster
 \enddefinition


\subsection{Algorithms}\label{s.algorithms}
In many arguments we use computer-aided test to check that a given graph satisfies certain conditions.
All
algorithms are
found
in \cite[Appendix]{Degtyarev.Rams}.
Given a graph
$\graph$, one can check whether:
\roster  \label{cond.tests}
\item\label{alg.prelim}
$\Fano(\graph)$ is hyperbolic, \ie, $\Gs_+\Fano(\graph)=1$; straight as it is, this
can often be done fast in bulk, see \cite[Lemmas~A.3, A.4]{Degtyarev.Rams};
\item\label{alg.master}
$\graph$ is extensible and admissible,
see
the master test in \cite[\S\,A.1.1]{Degtyarev.Rams};
\item\label{alg.geometric}
$\graph$ is
(sub-)geometric, see
\cite[\S\,A.1.2]{Degtyarev.Rams}.
\endroster

The principal tool of~\cite{Degtyarev.Rams} was starting from a sufficiently
large initial graph~$\graph_0$ and extending it by adding one vertex at a
time, see \S\,A.4 in \latin{loc.\ cit}. This is what is done in this paper,
too, see \autoref{S.quad} and \autoref{S.other}, where, without much
explanation, we merely state the updated results. The principal novelty of
this paper is \autoref{S.trig}, where, due to the more complicated geometric
nature of the problem, extra vertices have to be added in groups of up to
three. The details are discussed in \autoref{s.extending}.


\remark\label{rem.acceptable}
A technical, but crucial part of our (as well as any lattice-based)
approach is the fact that a geometric graph $\graph$ has $\rank\graph\le20$.
It follows that, if $\rank\graph'=20$, any geometric overgraph
$\graph\supset\graph'$ would be of the form
$\sat(\graph',\CK)$ for some $\CK \in \GK(\graph')$, and all such finite
index extensions of $\Fano(\graph')$ can easily be found
using~\cite{Nikulin:forms} (see the \emph{saturation lists} in
\cite[\S\,A.1.3]{Degtyarev.Rams}). Therefore, we introduce another technical
term: 
\[*
\text{$\graph$ is \emph{acceptable} if it is subgeometric and $\rank\graph\le19$}.
\]
It is understood that, after each step of every algorithm,
only
\emph{acceptable} graphs are left for the further processing (and, thus, we
do not need to add dozens of vertices to reach line counts over~$50$),
whereas each intermediate graph~$\graph$ of rank $\rank\graph=20$ is excluded upon
computing its saturation list and recording all ``interesting''
(\cf. \autoref{rem.trig.observed} below) geometric overgraphs to a global
master list.
\endremark



\subsection{Classification  of graphs in terms of girth}\label{s.classification}


Following~\cite{degt:lines,DIS}, we subdivide parabolic and hyperbolic
graphs~$\graph$ according to the
type of the minimal (in the sense of \autoref{conv.graphs}) affine Dynkin
diagram $\fiber\subset\graph$. The most important classes can also be
characterised in terms of
 the
 \emph{girth} $\girth(\graph)$
(the length of a
 shortest cycle in~$\graph$,
 with the convention that the girth of a forest is~$\infty$).
Thus, $\graph$ is called
 \roster*
 \item
 \emph{triangular}, or \smash{$\tA_2$}-, if $\girth(\graph)=3$,
 \item
 \emph{quadrangular}, or \smash{$\tA_3$}-,
 if $\girth(\graph)=4$,
 \item
 \emph{pentagonal}, or \smash{$\tA_4$}-,
 if $\girth(\graph)=5$,
 \item
 \emph{astral}, or \smash{$\tD_4$}-,
 if $\girth(\graph)\ge6$ and $\graph$ has a vertex
 of valency $\ge4$.
 \endroster
 All other graphs are \emph{locally elliptic},
 \ie, one has $\val v\le3$ for each
 vertex $v\in\graph$
 (and we assume $\girth(\graph)\ge6$ to exclude a few trivial
 cases).

Given a graph $\graph$
and a distinguished
connected parabolic subgraph $\fiber\subset\graph$,
the \emph{pencil} $\pencil$ induced by $\fiber$ is defined as
\[
\pencil:=\pencil(\graph\supset\fiber):=\fiber\cup
\bigl\{l\in\graph\bigm|\text{$l\cdot c=0$ for all vertices  $c\in\fiber$}\bigr\}.
\label{eq.pencil}
\]
This 
graph $\pencil\supset\fiber$ is parabolic as it is orthogonal to
$\kappa_\fiber$ (see \autoref{conv.graphs} and \autoref{rem.index}).
We have $\graph=\pencil\cup{\sec^*}$, where
\[
\sec^*:=\sec^*(\graph\supset\fiber):=
\bigl\{l\in\graph\sminus\fiber\bigm|
 \text{$l\cdot c=1$ for a vertex  $c\in\fiber$}\bigr\}.
\label{eq.sec*}
\]
The elements of $\sec^*$ are called the \emph{\rom(multi-\rom)sections}, or
\emph{$m$-sections} of~$\pencil$, where the integer $m:=l\cdot\kappa_\fiber$
is the \emph{multiplicity} of a section~$l$. If $m=1$, the section is called
\emph{simple}, otherwise, \emph{multiple}.
Fixing an order $\fiber=(c_1,\ldots,c_n)$, we also consider
\[
\alignedat2
\sec_i&:=\sec(\graph\supset\fiber\ni c_i)&&:=
\bigl\{l\in\graph\sminus\fiber\bigm|
 \text{$l\cdot c_j=\delta_{ij}$ for $c_j\in\fiber$}\bigr\},\\
\sec^*_i&:=\sec^*(\graph\supset\fiber\ni c_i)&&:=
\bigl\{l\in\graph\sminus\fiber\bigm|l\cdot c_i=1\bigr\}\supset\sec_i,
\endalignedat
\label{eq.sec}
\]
where $1\le i\le n$ and $\delta_{ij}$ is the Kronecker symbol.
We have
\[
\text{each set $\sec^*_i\supset\sec_i$, $i=1,\ldots,n$, is either elliptic or parabolic},
\label{eq.sec.elliptic}
\]
as it is
orthogonal to the isotropic vector $h-c_i\ne0$ (see \autoref{rem.index}).

In the future, we almost never use the correct,
but long notation referring to the full
flag, as $\graph\supset\fiber=(c_1,\ldots,c_n)$ are always assumed fixed.

\section{Triangular sets}\label{S.trig}

A \emph{triangular set}, or \emph{$\trig$-set}, is an induced subgraph of an
admissible graph whose all connected components are of type $\tA_2$,
$\bA_3$, $\bA_2$, or $\bA_1$. Clearly, each triangular pencil is a
$\trig$-set, but not \latin{vice versa}: a $\trig$-set may also be
elliptic, \ie, have no type~\smash{$\tA_2$} components. Combinatorially, a
$\trig$-set~$\set$
is uniquely of the form
\[
\ta_2\tA_2\splus a_3\bA_3\splus a_2\bA_2\splus a_1\bA_1,\qquad
(\ta_2,a_3,a_2,a_1)\in\NN^4,
\label{eq.coeff}
\]
and the \emph{coefficient quadruple} $(\ta_2,a_3,a_2,a_1)$ determines~$\set$ up
to isomorphism.
The isomorphism classes of $\trig$-sets (or coefficient
quadruples) are called \emph{patterns}.

We define the cardinality $\ls|\cs|$ of a pattern~$\cs$ as that of any of its
representatives and
introduce the following order on the set of patterns:
\[
\cs'\less \cs''\quad\text{iff}\quad\begin{cases}
 \ls|\cs'|<\ls|\cs''|\ \text{or}\\
 \ls|\cs'|=\ls|\cs''|\ \text{and}\ (\ta_2',a_3',a_2',a_1')>(\ta_2'',a_3'',a_2'',a_1''),
\end{cases}
\label{eq.order}
\]
where the coefficient quadruples are compared lexicographically. (Pay
attention to the \emph{reverse} lexicographic order!) This order is not to be
mixed with the perturbation order defined in \autoref{conv.graphs}.
The latter, for
$\trig$-sets, is easily described in terms of the coefficient quadruples:
$\cs'\pert\cs''$ if and only if $(\ta_2',a_3',a_2',a_1')$ is obtained from
$(\ta_2'',a_3'',a_2'',a_1'')$ by
a finite sequence of
\emph{elementary perturbations} of the form
\[*
\gathered
(\ta_2,a_3,a_2,a_1)\mapsto(\ta_2,a_3,a_2,a_1)+\delta,\\
\delta\in\bigl\{(0,0,-1,1),(0,-1,1,0),(0,-1,0,2),(-1,0,1,0)\bigr\}
\endgathered
\]
provided that, at each step, the quadruple remains in $\NN^4$.

\subsection{Constructing triangular graphs}\label{s.trig}
Let $\graph$ be a triangular admissible graph, \latin{cf.} \autoref{s.classification}.
Fix a type \smash{$\tA_2$} fiber
$\fiber=(c_1,c_2,c_3)\subset\graph$
and consider the pencil~$\pencil$ and sets $\sec_i$, $i=1,2,3$,
see \eqref{eq.pencil}
and~\eqref{eq.sec}, respectively.

\lemma\label{lem.3-section}
A triangular pencil~$\pencil$ has at most one multiple section\rom; hence,
\[*
0\le\ls|\graph|-\ls|\pencil\cup{\sec_1}\cup{\sec_2}\cup{\sec_3}|\le1.
\]
Furthermore, if a multiple section~$s$ exists, it is disjoint from each
$\sec_i$, $i=1,2,3$.
\endlemma

\proof
Consider a multiple section~$s$ and let $e:=h-c_1-c_2-c_3-s$.

If $s\cdot c_1=s\cdot c_2=s\cdot c_3=1$, then $e$ is orthogonal to~$h$ and
each $c_i$, $i=1,2,3$. Hence, $e=0$ (see \autoref{rem.index})
and any other line~$l$ intersects exactly
one of $c_1$, $c_2$, $c_3$, $s$, implying both statements.

Likewise, if $s\cdot c_1=s\cdot c_2=1$ and $s\cdot c_3=0$, then $e$ is an
exceptional divisor and $e\cdot s=e\cdot c_3=1>0$. Hence, any other line~$l$
intersects \emph{at most} one of $c_1$, $c_2$, $c_3$, $s$,
as otherwise $l\cdot e<0$ and $e$
would separate~$s$ and~$l$, see \eqref{eq.l34}.
\endproof

\lemma\label{lem.Delta}
Each set $\sec_i$, $i=1,2,3$, is a $\trig$-set.
\endlemma

\proof
In view of~\eqref{eq.sec.elliptic}, it suffices to rule out the connected
components of~$\sec_i$
containing $\tA_3$,
$\bA_4$, or $\bD_4$. Each of the offending graphs has a chain
\[*
\def\bul#1{\overset{\hbox to0pt{\hss$\scriptstyle s_#1$\hss}}{\mathrel\bullet}}
\def\edge{\joinrel\relbar\joinrel\relbar\joinrel}
\bul1\edge\bul2\edge\bul3
\]
and another vertex~$l$ adjacent to at least one of $s_1,s_2,s_3$.
Then $s:=s_1+s_2+s_3$ is a root, $s\cdot h=s\cdot c_i=3$, and $l\cdot s>0$;
hence, $e:=-h+c_i+s$ is an exceptional divisor separating $s_1$
and~$l$, see \eqref{eq.l34}.
\endproof

Now, our goal is to describe all geometric graphs $\graph$ of size
$\ls|\graph|\ge53$. In view of \autoref{lem.3-section}, it suffices to
consider trigonal pencils $\fiber\subset\pencil\subset\graph$ such that
\[*
\ls|\pencil|+\ls|\sec_1|+\ls|\sec_2|+\ls|\sec_3|\ge52.
\]
Clearly, we can also assume that $\pencil\subset\graph$ is maximal with
respect to~\eqref{eq.order} and the edges
$c_1,c_2,c_3$ of~$\fiber$ are ordered so that
${\sec_1}\grteq{\sec_2}\grteq{\sec_3}$. These assumptions give rise to the
following \emph{compatibility conditions}
on the patterns $\cp\ni\pencil$ and $\cs_i\ni\sec_i$:
\allowdisplaybreaks
\begin{align}
\label{eq.c1}
\cp\comp\cs_1 & \quad\text{if}\quad
 3\ls|\cs_1|+\ls|\cp|\ge52\ \text{and}\
 \bigl(\text{$\cs_1\lesseq\cp$ or $\cs_1$ is elliptic}\bigr);\\
\label{eq.c2}
(\cp\comp\cs_1)\comp\cs_2 & \quad\text{if}\quad
 2\ls|\cs_2|+\ls|\cs_1|+\ls|\cp|\ge52\ \text{and}\
 \cs_2\lesseq\cs_1;\\
\label{eq.c3}
(\cp\comp\cs_1\comp\cs_2)\comp\cs_3 & \quad\text{if}\quad
 \ls|\cs_3|+\ls|\cs_2|+\ls|\cs_1|+\ls|\cp|\ge52\ \text{and}\
 \cs_3\lesseq\cs_2.
\end{align}
In~\eqref{eq.c1}, we assume one of the terms, $\cp$ or~$\cs_1$, fixed
and treat the condition as a
restriction on
 the other term. The other two
conditions are
restrictions {\3 on}
the last term provided that the parenthesized
part is fixed.


Since $\trig$-sets appear as sets of sections, for an ordered fiber
$\tA_2\cong\fiber=(c_1,c_2,c_3)$ and $\trig$-set~$\set$ we
define
$\fiber\icup_i\set$, where $i=1,2,3$,
as the graph obtained from
  the disjoint union of $\fiber$ and $\set$
by connecting~$c_i$ to each vertex $v\in\set$
by a simple edge.

This construction extends to patterns, producing an isomorphism
class of graphs.
Checking the parameter quadruples $(\ta_2,a_3,a_2,a_1)$ one-by-one, it is
fairly easy to compute the sets of patterns
\[*
\aligned
\Pencils&:=\{\text{subgeometric triangular pencils}\}/{\cong},\rlap{ and}\\
\Sets&:=\{\text{$\trig$-sets~$\set$ such that \smash{$\tA_2\icup_i\set$} is subgeometric}\}/{\cong}.
\endaligned
\]
(To simplify the computation, for $\Pencils$ one can start from Shimada's
list~\cite{Shimada:elliptic} of Jacobian elliptic $K3$-surfaces,
and for
$\Sets$ one can take into account the bound $\val c_i\le20$ found
in~\cite{Veniani}.) Then, condition~\eqref{eq.c1} becomes a
 binary relation from~$\Pencils$ to $\Sets$.
The other conditions also descend to patterns, as do
the rank functions:
\[
\alignedat2
\rank(\cp)&=\ls|\pencil|-\ta_2(\pencil)+2,\quad&\pencil&\in\cp\in\Pencils\\
\rank(\smash{\tA_2}\icup_i\cs)&=\ls|\set|-\ta_2(\set)+4,&\set&\in\cs\in\Sets,
\endalignedat
\label{eq.rank}
\]
where $\ta_2(\cdot)$ stands for the number of parabolic
components.

The next lemma is an immediate consequence of this computation. For the last
statement, we merely list all geometric extensions
(\eg \, using algorithms from \cite[Appendix A]{Degtyarev.Rams})
of the four graphs
of rank~$18$ or~$19$, see~\eqref{eq.rank};
in fact, the sharp bound in \autoref{lem.Si}\iref{i.Si<=13} is $\ls|\graph|\le29$.

\lemma\label{lem.Si}
In a geometric graph~$\graph$, for any type~\smash{$\tA_2$} fiber~$\fiber$ one
has\rom:
\roster
\item\label{i.Si<=18}
$\ls|\sec_i|\le18$ for each $i=1,2,3$\rom;
\item\label{i.Si<=15}
if $\sec_i$ is elliptic, then $\ls|\sec_i|\le15$\rom;
\item\label{i.Si<=13}
if $\sec_i$ is elliptic and $\ls|\sec_i|>13$, then $\ls|\graph|\le52$.
\done
\endroster
\endlemma


\remark\label{rem.compatibility}
In view of \autoref{lem.Si}, the compatibility condition~\eqref{eq.c1}
simplifies to the form
\[*
\cp\comp\cs_1 \quad\text{if}\quad
 3\ls|\cs_1|+\ls|\cp|\ge52\ \text{and}\ \cs_1\lesseq\cp
\]
more consistent with~\eqref{eq.c2} and~\eqref{eq.c3}. Indeed, if $\cs_1$ is
elliptic, we can assume that $\ls|\cs_1|\le13$; then necessarily
$\ls|\cp|\ge13$, see~\eqref{eq.c1},
and $\cs_1\lesseq\cp$ holds automatically.
\endremark

\lemma\label{lem.excl-13}
If $\pencil\subset\graph$ is a maximal, with respect to~\eqref{eq.order},
triangular pencil and $\ls|\graph|>52$, then $\ls|\pencil|\ge14$.
\endlemma

\proof
As explained in \autoref{rem.compatibility}, in view of \autoref{lem.Si} we
have $\ls|\pencil|\ge13$ and
$\ls|\sec_3|\le\ls|\sec_2|\le\ls|\sec_1|\le\ls|\pencil|$.
If $\ls|\pencil|=13$, then $\ls|\sec_3|=\ls|\sec_2|=\ls|\sec_1|=13$
by~\eqref{eq.c1} and, moreover, $\pencil$ must have a multiple section~$s$.
Assuming that $s\cdot c_1=s\cdot c_2=1$, so that $c_1,c_2,s$ constitute a
triangle, the union $(c_1,c_2,s)\sqcup\sec_3\subset\graph$ is a triangular
pencil (see \autoref{lem.3-section}) with $16$ vertices,
contradicting the maximality of~$\pencil$.
\endproof

%
%
%

\subsection{\AAmple\ collections}\label{s.ample}
Let $\graph_0\supset\fiber=(c_1,c_2,c_3)\cong\smash{\tA_2}$ be a
triangular graph,
$i=\0,1,2,3$ 
 a parameter,
and $\cs\in\Pencils$ (if $i=\0$) or $\cs\in\Sets$ (otherwise) a pattern.
An overgraph $\graph\supset\graph_0$ is said to
\emph{represent} $(\graph_0,\cs)_i$ if
\[*
\alignedat2
&\text{$\sec(\graph\supset\fiber\ni c_i)\in\theta$ \ and \
 $\graph\sminus\sec_i=\graph_0$},\quad&&\text{if $i=1,2,3$}\\
&\text{$\pencil(\graph\supset\fiber)\in\cs$ \ and \
 $\graph\sminus\pencil=\graph_0\sminus\fiber$},&&\text{if $i=\0$}.
\endaligned
\]
A pair $(\graph_0,\cs)_i$ is called \emph{\aample}\ if it cannot be represented by an
acceptable graph, see \autoref{rem.acceptable}. Both notions extend to a
collection $\CG$ of graphs: $\graph$ represents $(\CG,\cs)_i$ if it
represents
$(\graph_0,\cs)_i$ for a graph $\graph_0\in\CG$, and $(\CG,\cs)_i$ is
\aample\ if so is each $(\graph_0,\cs)_i$, $\graph_0\in\CG$.
Iterating, we extend both notions to a compatible
 (\ie, satisfying the compatibility conditions from \autoref{s.trig})
 collection of patterns
\[*
\cp\in\Pencils,\quad \cs_i\in\Sets,\ i=1,\ldots,n\le3.
\]
Clearly, if $\cp\comp\ldots\comp\cs_n$ is
\aample, so is any $\cp'\comp\ldots\comp\cs'_n$ with
$\cp\pert\cp'$, \dots, $\cs_n\pert\cs'_n$.

Our proof of \autoref{thm-main}
essentially
reduces to applying the algorithm in \autoref{s.extending} below
to show
that any compatible
collection $\cp\comp\cs_1\comp\cs_2\comp\cs_3$ is \aample: indeed,
conditions~\eqref{eq.c1}--\eqref{eq.c3} guarantee that on the way
we will
encounter \emph{all} subgeometric graphs~$\graph$ such that either $\ls|\graph|>52$
or $\rank\graph=20$, and in the latter case it would suffice to analyze all
geometric saturations of~$\graph$. To this end, we introduce the inductive
notion of a \emph{ruled out collection}:
\roster*
\item
any \aample\ compatible collection is considered ruled out;
\item
in general, a compatible collection $\cp\comp\cs_1\comp\ldots\comp\cs_n$,
$n\le2$, is ruled out if so is \emph{any} compatible extension
$(\cp\comp\cs_1\comp\ldots\comp\cs_n)\comp\cs_{n+1}$, $\cs_{n+1}\in\Sets$.
\endroster
By a machine aided computation, we establish the following
statement;
its proof
is given in  \autoref{proof.trig}, after we collect all the necessary facts in \autoref{s.extending},
\autoref{s.multi-pat} and \autoref{s.aggressive}.

\proposition[see \autoref{proof.trig}]\label{prop.trig}
Each compatible pair $\cp\comp\cs_1$, where
\[*
\cp\in\Pencils_{14}:=\bigl\{\cp\in\Pencils\bigm|\ls|\cp|\ge14\bigr\}
 \quad\text{and}\quad
\cs\in\Sets,
\]
is ruled out.
\endproposition

By the very definition, the assertion of \autoref{prop.trig} means that,
	for each representative~$\graph$ of any compatible
	collection $\cp\comp\cs_1\comp\cs_2\comp\cs_3$, $\cp\in\Pencils_{14}$,
	one has
	either $\ls|\graph|=52$ or $\rank\graph=20$, and, moreover,
all such representatives are
	encountered in the course of the proof. The latter fact enables us to obtain
    the  complete list of representatives~$\graph$ of  compatible
	collections $\cp\comp\cs_1\comp\cs_2\comp\cs_3$, $\cp\in\Pencils_{14}$,
	such that
	$\ls|\graph|>52$ and $\rank\graph=20$ (see \autoref{add.list}).


\subsection{Extending a graph by a triangular set}\label{s.extending}
The heart of the computation is an algorithm extending a given subgeometric
graph $\graph_0$ by a given pattern $\cs\in\Sets$, the goal being listing
all subgeometric overgraphs $\graph\supset\graph_0$
representing $(\graph_0,\theta)_i$
(where $i=\0,1,2,3$ is also fixed, see \autoref{s.ample}).
The elements of
$\graph\sminus\graph_0$ are referred to as sections, whereas the
connected  components
 of $\graph\sminus\graph_0\in\cs$ are \emph{polysections}.
The algorithm is similar to that of~\cite{Degtyarev.Rams},
except that we
can no longer guarantee that the sections are pairwise disjoint.
Therefore, instead of adding to~$\graph$ one section at time, we fix~$\cs$ in
advance
and add whole polysections,
in the order \smash{$\tA_2$}, $\bA_3$, $\bA_2$, $\bA_1$.

\convention
From now on, following~\cite[\S\,A.3]{Degtyarev.Rams},
we identify a section~$v$ with its \emph{support}
\[*
\supp v:=\bigl\{u\in\graph_0\bigm|v\cdot u=1\bigr\}
\]
and thus treat it as a subset of~$\graph_0$. We also keep the
notation
\[*
\graph_0\sqcup\bv\quad\text{and}\quad\graph_0\sqcup\bv(\gm)
\]
for a multiset $\bv$ (which we no longer assume sorted) and
$(\ls|\bv|\times\ls|\bv|)$-matrix~$\gm$.
As explained in \autoref{s.algorithms}, only acceptable graphs are
retained after each step.
\endconvention


In practice, we start with computing the group $G_0:=\Aut\graph_0$
and set
\[
\CS(\graph_0):=\bigl\{s\subset\graph_0\bigm|s\cap\fiber=\fixed\bigr\}
\label{eq.CS}
\]
of sections of $\graph_0$ satisfying extra conditions imposed by the problem
at hand. (Here,
$\fixed\subset\fiber$ is a certain subset fixed in advance. We can
also take into account a few obvious geometric restrictions,
but this is not crucial: ``wrong'' sections are immediately ruled out by the
preliminary tests in \autoref{s.algorithms}\iref{alg.prelim}.
We
omit many other technical tweaks, referring to the code~\cite{degt.Rams.anc}
as the ultimate source.)
Then, running the tests
cited in \autoref{s.algorithms}\iref{alg.master} and~\iref{alg.geometric},
we compute the sets
\[
\aligned
\Sec_1(\graph_0)&:=\bigl\{v\in\CS(\graph_0)\bigm|
 \text{$\graph_0\sqcup v$ is acceptable}\bigr\},\\
\gm(\graph_0)&:=\bigl\{\bv\in\Sec_1(\graph_0)^{\dim\gm}\bigm|
 \text{$\graph_0\sqcup\bv(\gm)$ is acceptable}\bigr\},
\endaligned
\label{eq.S(graph)}
\]
where $\gm=\bA_2$, \smash{$\tA_2$}, $2\bA_1$, $\bA_3$
(in this order).
Certainly, the tests are applied to a single representative of each
$G_0$-orbit; in what follows (\cf., \eg, \autoref{rem.new.Sn})
this convention is taken for granted.
This computation is aborted
if a ``required'' list
is
empty (\eg, if
$\Sec_2(\graph_0)=\varnothing$ whereas $\cs$ contains
\smash{$\tA_2$} or $\bA_3$, \cf. the next remark).

\remark[a technical detail]\label{rem.A2}
The set $\Sec_2(\graph_0)$ is used in the computation of
\smash{$\tilde\Sec_2(\graph_0)$}: we consider only those triples
$(v_1,v_2,v_3)$ for which $(v_i,v_j)\in\Sec_2(\graph_0)$ for all
$1\le i<j\le3$. Likewise, both $\Sec_2(\graph_0)$ and $2\Sec_1(\graph_0)$ are
used in the computation of $\Sec_3(\graph_0)$. Furthermore,
$2\Sec_1(\graph_0)$ is used at all subsequent steps: when iterating
\[*
\graph_{n-1}\sqcup\bv(\gm_n)=
\bigl(\graph_0\sqcup\bu(\cdot)\bigr)\sqcup\bv(\gm_n),
\]
in~\eqref{eq.Sn} below,
we check first that $(u,v)\in2\Sec_1(\graph_0)$ for all $u\in\bu$,
$v\in\bv$.
\endremark

Now, let $\cs=\ta_2\smash{\tA_2}+a_3\bA_3+a_2\bA_2+a_1\bA_1$, so that
$N:=\ta_2+a_3+a_2+a_1$ is the number of components, and let
$\gm_1\ge\ldots\ge\gm_N$ be the types of the
components of~$\cs$
ordered \via\ $\smash{\tA_2}>\bA_3>\bA_2>\bA_1$).
Then,
we start from
$$\gS_0:=\{\graph_0\}$$ and run the computation in up to $N$ steps.

\subsubsection*{Step $n\ge1${\rm :}}
for each graph $\graph_{n-1}\in\gS_{n-1}$, we
pick a single representative~$\bv$ of each
$(\Aut\graph_{n-1})$-orbit on $\gm_n(\graph_0)$ and use the tests of
\autoref{s.algorithms}\iref{alg.master} and~\iref{alg.geometric}
to compute
\[
\gS_n(\graph_{n-1}):=\bigl\{\graph_\bv:=\graph_{n-1}\sqcup\bv(\gm_n)\bigm|
 \text{$\graph_\bv$ is acceptable}\bigr\}.
\label{eq.Sn}
\]
The step concludes by uniting all sets $\gS_n(\graph_{n-1})$,
$\graph_{n-1}\in\gS_{n-1}$, obtained followed by
retaining a single representative of
each graph isomorphism class.

The algorithm terminates  either upon the completion of Step~$N$ (resulting in a list
$\gS_N(\graph_0)$ to be processed by other means) or when one of the previous
steps results in an empty list $\gS_n(\graph_0)=\varnothing$,
implying that $(\graph_0,\cs)$
is \aample.

\remark\label{rem.digraph}
Both auto- and isomorphisms of graphs are computed using the \texttt{digraph}
package in \GAP~\cite{GAP4}. All morphisms are restricted: we assume the
fiber~$\fiber$ fixed as a set and, typically, one or two edges
$c_i\in\fiber$ fixed pointwise, so that the set $\CS(\graph_0)$ in
\eqref{eq.CS} be invariant.
\endremark

\remark[a technical detail]\label{rem.new.Sn}
Each time the matrix $\gm$ changes, \ie, whenever $\gm_n>\gm_{n+1}$, we
recompute the (relevant) sets $\gm(\graph_n)$ for each graph
$\graph_n\in\gS_n$ and use these new lists in the subsequent steps. Instead
of starting from scratch, as in the case of $\graph_0$, we merely run the
tests on the ready lists $\gm(\graph_m)$ for the last subgraph
$\graph_m\subset\graph_n$ for which they have been computed.
\endremark

\subsection{Processing several patterns}\label{s.multi-pat}
The material of
this section is of a purely technical nature; however,
it is the tweak described
here that makes the computation much faster and eventually helps it
to terminate reasonably fast.

Typically, we fix a subgeometric graph~$\graph_0$ and try to rule out a whole
collection of patterns $\Sets(\graph_0)$. Since patterns tend to have similar
initial sequences, processing them all one-by-one would force us to repeat the
same steps of the computation over and over again. To avoid the repetition
and remove a number of redundant steps,
we sort the patterns in the \emph{direct} lexicographic order and process
them simultaneously, organizing the computation into four layers: the
outermost~\smash{$\tA_2$}, $\bA_3$, $\bA_2$, and the innermost~$\bA_1$.

Each inner layer starts from a certain intermediate graph~$\graph$ and
processes a collection of patterns $\Sets(\graph)$. If the algorithm
terminates prematurely, at a certain pattern $\cs$, we conclude that
$(\graph,\cs)$ is \aample, and hence so is $(\graph,\cs')$
whenever $\cs\pert\cs'$. This fact is reported to the previous layer, where
the information is consolidated and often results in excluding the graph~$\graph$
and/or some patterns from the further consideration. We refer to the
code~\cite{degt.Rams.anc} for
the precise details (we implement each next layer as a hook within
the previous one, where it is used to modify the intermediate lists);
here, we merely illustrate the paradigm by the following simple example.

\example
Assume that the patterns to be considered are
\[*
\Sets(\graph_0)=\bigl\{
 2\bA_2\oplus6\bA_1,3\bA_2\oplus4\bA_1,4\bA_2\oplus2\bA_1\bigr\},
\]
so that only two layers
of computation are required. We run the first two
steps of the $\bA_2$-layer, resulting, say, in a list
$\gS_2=\{\graph_2',\graph_2''\}$, and switch to the $\bA_1$-layer
for each of the two graphs. Assume that this inner layer terminates at
\roster*
\item
step 4 for $\graph_2'$ $\Rightarrow$
 $(\graph_2',4\bA_1)$, $(\graph_2',\bA_2\oplus3\bA_1)$,
  $(\graph_2',2\bA_2\oplus2\bA_1)$ are \aample,
\item
step 5 for $\graph_2''$ $\Rightarrow$
 $(\graph_2'',5\bA_1)$, $(\graph_2'',\bA_2\oplus4\bA_1)$ are \aample.
\endroster
(Obviously, $4\bA_1\pert\bA_2\oplus3\bA_1\pert2\bA_2\oplus2\bA_1$ and
$5\bA_1\pert\bA_2\oplus4\bA_1$.) We conclude that $\graph_2'$ can be excluded
from $\gS_2$ and that both $2\bA_2\oplus6\bA_1$ and $3\bA_2\oplus4\bA_1$
can be excluded from~$\Sets(\graph_0)$. Therefore, we can run \emph{two} more
steps of the $\bA_2$-layer on the new reduced list $\{\graph_2''\}$, followed by the
$\bA_1$-layer on the result. (The $\bA_1$-layer after Step~$3$ can be skipped
as $\cs=3\bA_2\oplus4\bA_1$ has already been ruled out!)

If it were not for~$\graph_2''$ (\eg, if the $\bA_1$-layer terminated at a
step $\le4$ for each of the two graphs), we would have stopped immediately, as all
elements of $\Sets(\graph_0)$ would have been ruled out by the $\bA_1$-layer
after Step~2.
\endexample

\subsection{The aggressive version}\label{s.aggressive}
In certain cases,
one can argue that, in order to achieve the goal $\ls|\graph|\ge53$,
the overgraph
$\graph\supset\graph_0$ must be spanned over~$\graph_0$ by a few pairwise
disjoint vertices \emph{independent} over~$\graph_0$.
(Precisely, this condition means that $\Fano(\graph)\otimes\Q$ is
generated over $\Fano(\graph_0)\otimes\Q$ by $(\rank\graph-\rank\graph_0)$
pairwise disjoint vertices.) In this case, we switch to the
\emph{aggressive version} of the algorithm, \ie, we
\roster*
\item
add disjoint vertices only (the $\bA_1$-layer),
\item
disregard the extra vertices that do not increase rank, and
\item
check the saturation lists of all intermediate graphs, including~$\graph_0$,
\endroster
\cf. the \emph{progressive mode} in \cite[\S\,A.4.4]{Degtyarev.Rams}.

Precisely, this approach is used
\roster*
\item
in the proof of \autoref{lem.Si}\iref{i.Si<=13}: we add up to two disjoint
vertices, and
\item
at the steps $\fixed=\{c_2\}$ or $\{c_3\}$ in \autoref{proof.trig} below
(\cf. also \autoref{rem.c3}), when extending a graph~$\graph_0$ of the
submaximal rank $\rank\graph_0=19$.
\endroster

\subsection{Proof of \autoref{prop.trig}}\label{proof.trig}
As stated, the proof is an explicit machine aided computation
using the algorithm described in \autoref{s.extending} and
\autoref{s.multi-pat}. It runs in
two steps: first, for each pattern $\cs_1\in\Sets$, we rule out \emph{almost}
all compatible patterns $\cp\in\Pencils_{14}$; the set of
these patterns is denoted by $\Pencils^-_{14}(\cs_1)$.
Then, for each
$\cp\in\Pencils_{14}$, we rule out the remaining patterns $\cs_1\in\Sets$ such
that $\cp\comp\cs_1$ and $\cp\notin\Pencils^-_{14}(\cs_1)$.

\remark
The choice of the set $\Pencils^-_{14}(\cs_1)$ at the first step looks quite
arbitrary, and indeed so it is. \emph{As a rule}, we let
$\cp\in\Pencils^-_{14}(\cs_1)$ if $\cp\comp\cs_1$ and
\[*
\rank\cp\le\rank(\fiber\icup_{1}\cs_1).
\]
However, a few border cases are subject to further manual tweaking,
which is based on
experiments.
More precisely, depending on the values $r:=\rank(\fiber\icup_{1}\cs_1)$,
$n:=\ls|\cs_1|$, the following coefficient quadruples
$(\ta_2,a_3,a_2,a_1)$
are \emph{excluded} from $\Pencils^-_{14}(\cs_1)$:
\def\-{\:\quad}
\begin{align*}
(r,n)=
(12, 11)\- & (*, *, *, *), \cr
(12, 12)\- & (*, *, *, *), \cr
(13, 11)\- & (*, *, *, *), \cr
(13, 12)\- & (3, *, *, *), (2, *, *, *), (1, *, *, *), \cr
(13, 13)\- & (2, *, *, *), (1, *, *, *), \cr
(14, 11)\- & (*, *, *, *), \cr
(14, 12)\- & (3, *, *, *), (2, *, *, *), (1, *, *, *), \cr
(14, 13)\- & (2, *, *, *), (1, *, *, *), \cr
(14, 14)\- & (2, *, *, *), (1, *, *, *), \cr
(15, 11)\- & (*, *, *, *), \cr
(15, 13)\- & (2, 0, 0, *), (1, 0, *, *), \cr
(15, 14)\- & (2, 1, 0, *), (2, 0, *, *), (1, 0, *, *), \cr
(15, 15)\- & (1, *, *, *)
\end{align*}
(where, as usual, $*$ means any value).
\endremark

At the first step, we start with the graph
\[*
\graph_0:=\fiber\icup_1\set_1,\quad \set_1\in\cs_1\in\Sets,
\]
and use \autoref{s.multi-pat} to find all graphs $\graph$ representing
$(\graph_0,\cp)_\0$, $\pi\in\Pencils^-_{14}(\cs_1)$.
Technically, we extend $\graph_0$ by a pattern~$\cs$ such that
\smash{$\tA_2\oplus\cs=\cp$}.
Thus, we let $\fixed=\varnothing$ in~\eqref{eq.CS} and use graph
auto-/isomorphisms
preserving $\fiber\ni c_1$ (see \autoref{rem.digraph}).
For each graph $\graph$ on the resulting list $\gS_N$, we consider the full
set $\Sets(\graph)$ of
patterns $\cs_2\in\Sets$ satisfying
\eqref{eq.c2}, and run the same algorithm with $\fixed=\{c_2\}$ and graph
morphisms preserving $c_1$ and $c_2$.

\remark\label{rem.no.c3}
Strictly speaking, we should have run the algorithm once more, using
$\fixed=\{c_3\}$ and patterns $\cs_3$ satisfying~\eqref{eq.c3}.
However, our thresholds are chosen so that the list $\gS_N$
resulting from the first run consists of relatively few graphs of rank~$19$
(for which the aggressive version is used, see \autoref{s.aggressive}) and very
few graphs of rank~$18$, for which the algorithm terminates fast and
rules everything out.
\endremark

At the second step, we start with a pencil
\[*
\graph_0:=\pencil\in\cp\in\Pencils_{14}
\]
and use \autoref{s.multi-pat} to find all graphs $\graph$ representing
$(\graph_0,\cs_1)_1$, $\cp\comp\cs_1$ and $\cp\notin\Pencils^-_{14}(\cs_1)$.
We let $\fixed=\{c_1\}$ in~\eqref{eq.CS} and use graph
morphisms preserving $\fiber\ni c_1$. As above, the relatively few graphs
obtained, all of rank~$19$ or~$18$, are ruled out by the
next run, using $\fixed=\{c_2\}$ and $\cs_2\in\Sets$ compatible in the sense
of~\eqref{eq.c2}.
This completes the proof of \autoref{prop.trig},
as well as of \autoref{add.list}
below.
\qed

\medskip
As explained
right after \autoref{prop.trig},
before discarding a graph~$\graph$ of
rank~$20$, we analyze its geometric finite index extensions and record those
of size greater than~$52$. The result of this analysis is stated below.

\addendum\label{add.list}
Let
$\graph$ be a geometric representative of a compatible collection
$\cp\comp\cs_1\comp\cs_2\comp\cs_3$, $\cp\in\Pencils_{14}$, such that $\rank\graph=20$
and $\ls|\graph|>52$. Then $\graph$ is one of the eight smooth configurations
found in~\cite{DIS}, see the first eight rows of \autoref{tab.smooth}.
\done
\endaddendum

\remark\label{rem.trig.observed}
In
order to produce a plethora of examples of large configurations of lines,
when discarding the graphs of rank~$20$ (see \autoref{rem.acceptable}) we
collected all extended graphs with at least $48$ lines or at least six
exceptional divisors; the results are found in \cite{degt.Rams.anc}. In particular, in
addition to the surfaces listed in \autoref{tab.smooth},
we found but one quartic with $52$ lines
and two nodes and two quartics with $50$ lines and one node each. Besides,
there are quite a few quartics with non-empty singular locus and $48$ lines,
suggesting once again that $48$ is a reasonable threshold to cut the
classification.
\endremark


\section{Smooth quartics}\label{S.smooth}

In
this section we temporarily assume that the polarized
lattice  $S \ni h$ contains no exceptional divisors;
this will be used later
in our study of
line configurations on smooth quartic surfaces
$X_4\subset\Cp3$.

In this case, we need to change the notion of admissible lattice/graph.
Namely, a polarized lattice $S\ni h$ is called \emph{smooth}, or
\emph{$s$-admissible}, if it contains neither exceptional divisors ($e^2=-2$,
$e\cdot h=0$) nor $2$-isotropic vectors ($e^2=0$, $e\cdot h=2$). A
graph~$\graph$ is \emph{smooth}, or \emph{$s$-admissible}, if so is the
lattice $\Fano(\graph)$.

If $S\ni h$ is smooth, then $\rt(S,h)=\varnothing$; hence, automatically,
$\Weyl=\varnothing$ in \eqref{eq-def-Fano} and we have a well-defined Fano
graph
\[
\Fn(S,h):=\Fn_\varnothing(S,h)=\root_1(S,h).
\label{eq.Fn.smooth}
\]
An crucial consequence of~\eqref{eq.Fn.smooth} is the fact that $\Fn$ is
monotonous:
\[
\text{if $S'\supset S\ni h$ are smooth, then $\Fn(S',h)\supset\Fn(S,h)$}.
\label{eq.monotonous}
\]

\lemma\label{lem.smooth}
Let $\graph$ be a smooth graph and $\fiber=(c_1,c_2,c_3)\subset\graph$ a
triangle. Then\rom:
\roster
\item\label{i.3-section}
$\fiber$ has a unique $3$-section $c_0:=h-c_1-c_2-c_3\in\graph$\rom;
\item\label{i.components}
the pencil $\pencil$ as in~\eqref{eq.pencil} and $\trig$-sets $\sec_i$ as
in~\eqref{eq.sec} have no connected components of types $\bA_2$ or $\bA_3$.
\endroster
\endlemma

\proof
Clearly, $c_0$ as in Statement~\iref{i.3-section} is a $3$-section in the
lattice spanned by $h$ and~$\fiber$. By~\eqref{eq.monotonous}, it remains a
$3$-section in any larger \emph{smooth} lattice/graph.

If the pencil~$\pencil$ has a pair $(s_1,s_2)\cong\bA_2$,
by~\eqref{eq.monotonous} it also has
$s_3:=\kappa_\fiber-s_1-s_2$, so that $(s_1,s_2,s_3)\cong\smash{\tA_2}$.
If $\pencil$ has $(s_1,s_2,s_3)\cong\bA_3$, then $\kappa_\fiber-s_1-s_2-s_3$
is an exceptional divisor. The same argument applies to each set $\sec_i$,
$i=1,2,3$, except that we replace $\kappa_\fiber$ with
$(c_0+\ldots+c_3)-c_i=\kappa_\fiber+c_0-c_i$.
\endproof

In view of \autoref{lem.smooth}\iref{i.3-section}, each type~$\tA_2$ fiber
$\fiber_0:=\fiber=(c_1,c_2,c_3)$ gives rise to three more, \viz.
$\fiber_i:=(c_0,\ldots,\hat c_i,\ldots)$, $i=1,2,3$ (where, as usual,
$\hat c_i$ indicates that $c_i$ has been omitted). Thus, we can shift the
paradigm and, instead of considering a pencil~$\pencil$ and three sets
$\sec_i$, we can speak about ``blending'' four pencils
\[*
\pencil_0:=\pencil=\pencil(\graph\supset\fiber_0),\qquad
\pencil_i:=\fiber_i\sqcup\sec_i=\pencil(\graph\supset\fiber_i),\quad
 i=1,2,3.
\]
Assuming, as above, that $\pencil$ is a maximal pencil in~$\graph$, we can
replace \eqref{eq.c1}--\eqref{eq.c3} with a stronger set of compatibility
conditions: 
\begin{align}
\label{eq.s1}
\cp\comp\cs_1 & \quad\text{if}\quad
 3\ls|\cs_1|+\ls|\cp|\ge48\ \text{and}\
 \cs'_1\lesseq\cp;\\
\label{eq.s2}
(\cp\comp\cs_1)\comp\cs_2 & \quad\text{if}\quad
 2\ls|\cs_2|+\ls|\cs_1|+\ls|\cp|\ge48\ \text{and}\
 \cs'_2\lesseq\cs'_1;\\
\label{eq.s3}
(\cp\comp\cs_1\comp\cs_2)\comp\cs_3 & \quad\text{if}\quad
 \ls|\cs_3|+\ls|\cs_2|+\ls|\cs_1|+\ls|\cp|\ge48\ \text{and}\
 \cs'_3\lesseq\cs'_2.
\end{align}
Here, $'$ stands for the operator $\set\mapsto\set':=\smash{\tA_2}\sqcup\set$
and its natural descent to the set of patterns. (Recall also that, for
\autoref{th.smooth}, we need to change the threshold to $\ls|\graph|\ge49$.)
In particular,
from~\eqref{eq.s1}--\eqref{eq.s3} we immediately
conclude that 
\[ \label{eq-twelve}
\ls|\pencil|\ge 15,
\]
\cf. \autoref{lem.excl-13}.
Indeed, taking into account the $3$-section given by \autoref{lem.smooth},
we can rewrite~\eqref{eq.s3} in the form
$\ls|\cs'_3|+\ls|\cs'_2|+\ls|\cs'_1|+\ls|\cp|\ge57=48+9$, and it remains to
observe that
the assumption $\cs'_3\lesseq\cs'_2\lesseq\cs'_1\lesseq\cp$
implies $\ls|\cs'_3|\le\ls|\cs'_2|\le\ls|\cs'_1|\le\ls|\cp|$.
Thus, unlike \autoref{prop.trig}, we do not need to introduce an analogue
$\Pencils_{15}\smooth$ of the set $\Pencils_{14}$:
the lower bound~\eqref{eq-twelve} would follow from the
compatibility assumptions.

Now, a computation similar to (but much faster than) that of \autoref{S.trig}
yields the following result (\cf. \autoref{prop.trig};
we retain the terminology of \autoref{S.trig}
and denote by $\Pencils\smooth\subset\Pencils$ and
$\Sets\smooth\subset\Sets$ the sets 
of patterns appearing in smooth graphs).

\proposition[\latin{cf.} \autoref{proof.trig}]\label{prop.smooth}
Each compatible  collection
$\cp\comp\cs_1\comp\cs_2\comp\cs_3$, as in \eqref{eq.s1}--\eqref{eq.s3},
with $\cp\in\Pencils\smooth$
and
$\cs_1, \cs_2, \cs_3 \in \Sets\smooth$,
either is ruled out or
has
 one of the five graphs
\[ \label{eq-monstrous-five}
\bZ_{49},\quad  (50),\quad  \bZ_{50},\quad (52),\quad   \bZ_{52}
\]
\rom(see \autoref{tab.smooth}\rom) as a geometric representative.
\done
\endproposition

Similar 
to \autoref{prop.trig}, this statement means that, with the five
exceptions listed in~\eqref{eq-monstrous-five},
any geometric representative~$\graph$ of a collection as in the
hypotheses has $\rank\graph=20$ and, moreover,
all such representatives
of rank~$20$ are encountered
and discarded
in the course of the proof.
As in \autoref{S.trig}, prior to discarding a graph
we compute all its
geometric finite index extensions, thus arriving at
the following
complete list
of geometric rank~$20$ Fano graphs $\graph$ with $\ls|\graph|\ge48+1$
(the extra~$1$ standing for the $3$-section given by \autoref{lem.smooth}).

%

\addendum\label{add.list.smooth}
Let
$\graph$ be a geometric representative of a compatible collection
$\cp\comp\cs_1\comp\cs_2\comp\cs_3$ as in \eqref{eq.s1}--\eqref{eq.s3},
where $\cp\in\Pencils\smooth$,
$\cs_1, \cs_2, \cs_3 \in \Sets\smooth$,
and $\rank\graph=20$.
Then $\graph$ is one of the $21$ smooth
rank~$20$
configurations
found in~ \autoref{tab.smooth}.
\done
\endaddendum

\remark\label{rem.c3}
Due to the lower threshold $\ls|\graph|\ge49$, occasionally we do have to run
the algorithm till the very last step $\fixed=\{c_3\}$,
\cf. \autoref{rem.no.c3}. This last step is quite expensive (as the set of
sections to begin with is quite large), but
fortunately it has to be done for four graphs only. This is
yet another indication of the fact that taking the classification down to
$48$ or fewer lines is hardly feasible.
\endremark

\remark\label{rem.smooth.full}
In the smooth case, we can further reduce the overcounting by a number of
tricks based on the monotonicity property~\eqref{eq.monotonous},
using all lines present in the
lattice rather than only those added explicitly.
Most notably,
before switching to a next step $\fixed=\{c_i\}$, $i\le3$, we
can replace each graph~$\graph$ with the union
\[*
\pencil\cup\sec_1\cup\ldots\cup\sec_{i-1}\quad\text{in}\quad
 \Fn(\Fano(\graph))
\]
and then retain a single representative of each isomorphism class of the list
obtained. Furthermore, in the subsequent computation we can impose an extra
condition that the above union should remain fixed.
We refer
to~\cite{degt:lines} for further details.
\endremark

\section{Quadrangular graphs}\label{S.quad}

In this section the \KK-quartic $X_4$ is allowed to have singular points as in \autoref{S.trig}.
Consider an $\tA_3$-graph~$\graph$
and fix a quadrangle $\fiber:=(c_1,c_2,c_3,c_4)\subset\graph$.
We assume the
edges $c_i\in\fiber$ numbered cyclically and ordered consecutively:
\[*
c_{i+4}=c_i,\qquad
c_i\cdot c_{i\pm1}=1,\quad c_i\cdot c_{i\pm2}=0.
\]
In addition to~\eqref{eq.pencil}--\eqref{eq.sec}, we consider the sets
\[*
\sec^*_{ij}:={\sec^*_i}\cup{\sec^*_j},\quad i=j\pm2.
\]

The assumption that $\girth(\graph)=4$ implies that
\[*
{\sec^*_{13}}\cap{\sec^*_{24}}=\varnothing\qquad\text{and}\qquad
\text{each graph $\sec^*_i$ is discrete}.
\]
The elements of each $\sec_i$ are
simple sections of~$\pencil$, and those of
\[*
{\sec^*_i}\sminus{\sec_i}={\sec^*_j}\sminus{\sec_j}
 ={\sec^*_i}\cap{\sec^*_j},\quad j=i\pm2,
\]
are double
sections; the set of all (multi-)sections of $\pencil$ is
$\sec^*={\sec^*_{13}}\cup{\sec^*_{24}}$.

From now on, we state our results for \emph{geometric} rather than subgeometric
graphs; in other words, we assume that $\graph$ admits a \emph{triangle free}
geometric saturation.

\lemma\label{lem.quad.sections}
In a geometric quadrangular graph $\graph\supset\fiber$ as above, one has
\roster
\item\label{i.Ei}
$\ls|\sec_i|\le10$ and $\ls|\sec^*_i|\le10$ \rom(sharp bounds\rom),
\item\label{i.Eij}
$\ls|\sec^*_{ij}|\le20$ \rom(a sharp bound\rom), and
\item\label{i.quad.sections}
$\ls|\sec^*|\le32$
\rom(the best known example being~$30$\rom).
\endroster
\endlemma

\proof
Statement~\iref{i.Ei} is a direct computation: since each $\sec^*_i$ is
discrete, the union $\fiber\cup{\sec^*_i}$ depends on just two parameters,
which can take the following values:
\[
\hbox{\hss$\vcenter{\halign{\hss$#$&$#$\hss&$#$
 &&\hbox to16pt{\hss$#$\hss}\cr
p_i&{}:=\ls|{\sec_i}|&{}=\!&
 0& 1& 2& 3& 4& 5& 6& 7& 8& 9&10\cr
q_i&{}:=\ls|{\sec^*_i}\sminus\set_i|&{}\le\!&
 8& 7& 6& 5& 5& 4& 4& 3& 2& 1& 0\cr
\crcr}}$\hss}
\label{eq.quad.table}
\]
The bound in~\iref{i.Eij} follows from~\iref{i.Ei} (assuming
$p_i\ge p_j$, we have $\ls|\sec^*_{ij}|\le2p_i+q_i$),
and its sharpness is established by an explicit construction (obtained in the
next computation). Finally, statement~\iref{i.quad.sections} is also obtained
by a computation similar to \cite[\S\,B.5]{Degtyarev.Rams}: assuming that
\[
\ls|\sec_1|\ge\ls|\sec_3|,\quad
\ls|\sec_2|\ge\ls|\sec_4|,\quad
\ls|\sec^*_{13}|\ge\ls|\sec^*_{24}|,
\label{eq.quad.compat}
\]
we start with a ``standard'' graph
$\fiber\cup{\sec^*_1}$, letting $(p_1,q_1)=(10,0)$, $(9,1)$, $(9,0)$, $(8,2)$,
$(8,1)$, or $(7,3)$, and build all possible consecutive extensions
\[
\graph:=(\fiber\cup{\sec^*_1})\cup{\sec_3}\cup{\sec_2}\cup
 ({\sec^*_2}\sminus{\sec_2})\cup{\sec_4}
\label{eq.quad.steps}
\]
\via\ discrete sets;
the minimal size of each set to be added is determined
using \eqref{eq.quad.table}, \eqref{eq.quad.compat}, and the goal
$\ls|\graph|\ge37$.
In most cases, this algorithm terminates (meaning that each sufficiently large
graph admitting a geometric \emph{triangle free} saturation is unacceptable, \cf.
\autoref{rem.acceptable}) at the very first nontrivial step~$\sec_3$; in very
few cases we also need to use $\sec_2$.
\endproof

\remark
In the proof of \autoref{lem.quad.sections}
and \autoref{prop.quad} below, the essential difference
from~\cite{Degtyarev.Rams} is that we do not limit the number of lines
intersecting both given ones (the presence of biquadrangles and longer
``polyquadrangles''); this makes the computation slightly more involved.
\endremark

\proposition\label{prop.quad}
For a geometric quadrangular graph~$\graph$ one has $\ls|\graph|\le48$.
\endproposition

\remark
It is unlikely that the bound given by \autoref{prop.quad} is sharp: the best
example that we found has $39$ lines
(see also \cite[Proposition 6.14]{Degtyarev.Rams}).
\endremark

\proof[Proof of \autoref{prop.quad}]
The proof is a computation similar to \cite[\S\,C.3]{Degtyarev.Rams}.
In view of \autoref{lem.quad.sections}, it suffices to consider all
quadrangular pencils $\pencil$ of size $\ls|\pencil|\ge17$; the list of such
pencils admitting a polarization is compiled
using~\cite{Shimada:elliptic}. Then, under the assumptions
of~\eqref{eq.quad.compat}, we must have $4p_1+2q_1+\ls|\pencil|\ge49$.
Similar to~\eqref{eq.quad.steps}, we start from a pencil~$\pencil$ and build
a list of consecutive acceptable extensions
\[*
\graph:=\pencil\cup{\sec_1}\cup({\sec^*_1}\sminus{\sec_1})\cup{\sec_3}\cup
 {\sec_2}\cup({\sec^*_2}\sminus{\sec_2})\cup{\sec_4}.
\]
Due to our modest goal $\ls|\graph|\ge49$,
\latin{de facto} the algorithm terminates at the first or, occasionally,
second step, so that we never need to consider even $\sec_3$.
\endproof


\section{Other types of graphs} \label{S.other}

For the other types (in the sense of \autoref{s.classification}) of graph,
the computation runs exactly as in~\cite{Degtyarev.Rams}, and we merely state
the updated results below. Remarkably, the upper bounds obtained are exactly
the same as in the smooth case (see~\cite{degt:lines}); furthermore, unlike
the case of octics (see~\cite{Degtyarev.Rams}), all extremal configurations
are smooth.

\subsection{Pentagonal graphs}\label{s.pent}
Recall that the assumption $\girth(\graph)=5$ implies that, for a fiber
$\smash{\tA_4}\cong\fiber\subset\graph$, one has
\roster*
\item
$\sec^*={\sec_1}\cup\ldots\cup{\sec_5}$,
\ie, all sections are simple;
\item
each graph $\sec_i^*=\sec_i$, $i=1,\ldots,5$, is discrete.
\endroster

The following bounds are sharp, and there are but two
geometric pentagonal
graphs with $30$ vertices, \viz. $\Phi'_{30}$ and $\Phi''_{30}$ (see
\cite{degt:lines}). Both represent configurations of lines on \emph{smooth}
quartic surfaces only.

\lemma[\cf. {\cite[Lemma 7.6]{Degtyarev.Rams}}]\label{lem.pent}
Let $\graph$ be a geometric pentagonal graph, and let $\fiber\subset\graph$
be a
type~\smash{$\tA_4$} subgraph. Then\rom:
\roster
\item
one has $\ls|\sec_*|\le16$\rom;
\item
if $\ls|\sec_*|\ge14$, then $\ls|\graph|\le29$.
\done
\endroster
\endlemma

\proposition\label{prop.pentagonal}
One has $\ls|\graph|\le30$ for any geometric pentagonal graph~$\pencil$.
\endproposition

\proof
In view of \autoref{lem.pent},
it suffices to consider pentagonal pencils~$\pencil$
such that $\ls|\graph|\ge17$;
they can be found using~\cite{Shimada:elliptic}.
\endproof

\subsection{Astral graphs}\label{s.astral}
We number the vertices $(c_1,\ldots,c_5)$ of a type $\bD_4$ fiber~$\fiber$
so that the \emph{central vertex}~$c_1$ is the one of valency~$4$
in~$\fiber$.
Recall that the assumption $\girth(\graph)\ge6$ implies that, for a fiber
$\smash{\tD_4}\cong\fiber\subset\graph$ one has
\roster*
\item
$\sec^*={\sec_1}\cup\ldots\cup{\sec_5}$, and
\item
the graph $\sec^*$ is discrete.
\endroster
Note though that it is \emph{not} true that all sections are simple: the
elements of $\sec_1$ are double sections.
(Recall that $\kappa_\fiber=2c_1+c_2+\ldots+c_5$.)

The following bounds are sharp, and the only geometric astral graph
with $27$ vertices
is $\Delta'_{27}$ (see \cite{degt:lines}); it is represented by a unique
\emph{smooth} quartic surface.

\lemma[\cf. {\cite[Lemma 7.3]{Degtyarev.Rams}}]\label{lem.astral}
Let $\graph$ be a geometric astral graph and $\fiber\subset\graph$ a
type~\smash{$\tD_4$} subgraph whose central vertex~$c_1$ has maximal valency
in~$\graph$. Then\rom:
\roster
\item
one has $\ls|\sec^*|\le12$\rom;
\item
if $\ls|\sec^*|\ge11$, then $\ls|\graph|\le27$.
\done
\endroster
\endlemma

\proposition\label{prop.astral}
One has $\ls|\graph|\le27$ for any geometric astral graph~$\pencil$.
\endproposition

\proof
In view of \autoref{lem.astral},
it suffices to consider astral pencils~$\pencil$
with $\ls|\graph|\ge18$;
they can be found using~\cite{Shimada:elliptic}.
\endproof

\subsection{Locally elliptic graphs}\label{s.elliptic}
Let us recall, that the case of locally elliptic graphs was considered in \cite{Degtyarev.Rams}. We have the
inequality
\[ \label{eq-loc-elliptic}
\ls|\graph|\le29
\]
for all geometric locally elliptic graphs $\graph$
(see \cite[(7.1)]{Degtyarev.Rams}).
Machine-aided experiments
suggest that the sharp bound is $\ls|\graph|\le25$ with a unique graph $\Lambda_{25}$
that attains the maximum , but such considerations are of no importance for the proof of \autoref{thm-main}.




\section{Proofs}\label{S.Proof}

In order to render our exposition self-contained, we recall
certain results from \cite{Degtyarev.Rams} in \autoref{s.linesonquartics},
before presenting the proofs
of the main results of the paper in \autoref{proof.main}, \autoref{proof.smooth}.

\subsection{Fano graphs of  \KK-quartics} \label{s.linesonquartics}

Let
$X_4 \subset \Cp{3}$ be a complex degree-$4$ surface
with at worst
Du Val (\latin{aka} $\bA$--$\bD$--$\bE$, or simple) singularities and let
$\pi\:\resxf \rightarrow X_4$ be the minimal resolution of its singularities.
Denote $\hp := \hhp\in\NS(\resxf)$.

Recall that $\resxf$ is a $K3$-surface. In particular,
given an irreducible curve  $C \subset \resxf$, we have
\[
\text{$C\cdot\hp = 1$ and $C^2 = -2$ if and only if
	$\pi(C)$ is a degree one curve on $X_{4}$}.
\label{eq-lines-quartic}
\]
The curves that satisfy  \eqref{eq-lines-quartic} are called
\emph{lines} on $\resxf$; they are obviously smooth and rational.
We follow \cite{degt:lines}
and define the \emph{\rom(plain\rom) Fano graph} of the quartic
$X_4$ as the loop free graph with vertices
\[ \label{eq-fanograph-def}
\graphl(X_4) := \bigl\{
\text{$(-2)$-curves $C\subset \resxf$  with $C\cdot\hp = 1$}\bigr\}
\]
and each pair of vertices $v, w\in  \graphl(X_4)$  connected by
an edge of multiplicity $v \cdot w$.
(Here and below, we always consider the intersection form "$\cdot$" on
$\NS(\resxf)$.)

Recall that, by \cite[(4.5)]{Degtyarev.Rams},
\begin{equation} \label{eq-fano-must-be-hyperbolic}
\mbox{the graph $\graphl(X_4)$ of a \KK-quartic with at least 25 lines is hyperbolic.}
\end{equation}

General theory of  lattice-polarized \KK-surfaces
(Nikulin~\cite{Nikulin:finite.groups},
Saint-Donat~\cite{Saint-Donat}; \cf. also
\cite[\thmcit~3.11]{DIS}
and \cite[\thmcit~7.3]{degt:singular.K3})
yields the following
statement.
(As in
\cite[Convention~1.4]{Degtyarev.Rams},
we say that 
the lattice $\NS(\resxf)$ is
\emph{\qgen by lines}
if it is
a finite index extension
of its sublattice generated by the classes of
lines on $\resxf$
\emph{and the quasi-polarization~$\hp$},
\ie, it is spanned by lines and~$h$ over~$\Q$.)

\theorem[see {\cite[Theorem 3.9]{Degtyarev.Rams}}]\label{th.K3quartics}
A graph $\graph$ is geometric
if and only if one has
$\graph\cong\graphl(X_4)$ for a quartic~$X_4$ such that
$\NS(\resxf)$ is \qgen by lines.
\done
\endtheorem


\remark \label{rem-lines-gen-enough}
As in the case of \KK-octics, by \cite[Lemma~2.8]{Degtyarev.Rams},
in order to study
the
\emph{maximal number of lines} we can restrict our attention to the case when  the lattice $\NS(\resxf)$ is
\qgen by lines and exceptional divisors (see also  \cite[\S\,8]{Degtyarev.Rams})
\endremark

Obviously, for a quartic $X_4$  with non-empty singular locus
the graph $\graphl(X_4)$ does not completely describe
the configuration of lines on $X_4$.
The latter can be inferred from
the
bi-colored \emph{extended Fano graph}
\[  \label{eq-def-exFano-surface}
\graphe(X_4) := \bigl\{
\text{$(-2)$-curves $C\subset \resxf$  with $C\cdot\hp \le 1$}\bigr\},
\]
with the colour of each vertex $C$ defined as $C\cdot h$.
For such graphs we have a more general statement.

\theorem[see {\cite[Theorem 3.10]{Degtyarev.Rams}}]\label{th.K3quartics.ext}
A bi-colored graph~$\graph'$
is geometric
if and only if $\graph'\cong\graphe(X_4)$ for a
\KK-quartic $X_4 \subset \Cp{3}$.
\done
\endtheorem

\subsection{Proof of \autoref{thm-main}}\label{proof.main}

By \autoref{th.K3quartics.ext}, the
assertion of \autoref{thm-main} is equivalent to the  statement
that there are no
\begin{equation} \label{eq-graph-main}
\mbox{geometric bi-colored graphs
$\graph'$ such with $\ls|\spp_1\graph'|\geq53$ and
 $\ls|\spp_0\graph'|\ne0$,}
\end{equation}
where $\spp_j\graph'$ stands for the induced subgraph of $\graph'$
given by all its vertices of color $j = 0,1$.
Moreover, by
\eqref{eq-fano-must-be-hyperbolic},
the (plain) graph $\graph:=\spp_1\graph'$ is  a
$\fiber$-graph for a
certain affine Dynkin diagram \smash{$\fiber \ge \tA_2$}.

The monotonicity given by
\cite[Lemma~2.8]{Degtyarev.Rams} combined with the consideration of
\autoref{S.other} implies that the graph
$\graph$ is neither pentagonal
(see \autoref{prop.pentagonal}),
nor astral (see \autoref{prop.astral}), nor locally elliptic
(see \eqref{eq-loc-elliptic}). Finally, \autoref{prop.quad} shows that
$$
\mbox{$\graph$ is triangular. }
$$
Then,
by \autoref{lem.excl-13} and \autoref{prop.trig},
we
necessarily
have
$\rank(\graph) = 20$,
upon which \autoref{add.list} implies that $\graph$ is one of the eight
smooth configurations found in~\cite{DIS}.
For each
graph~$\graph$ obtained
we compute its  extended saturation(s) $\satex\graph$
and check that none of the resulting bi-colored graphs
has a vertex of color zero
(exceptional divisor),
completing
the proof of
 \autoref{thm-main}.  \hfill $\Box$

\subsection{Proof of \autoref{th.smooth}}\label{proof.smooth}
When dealing
with
smooth quartics,
we can confine ourselves
to the case where the lattice $\NS(X_4)$ is
\qgen by lines,
see \autoref{rem-lines-gen-enough} and
\eqref{eq.monotonous}.
\autoref{th.K3quartics} reduces the proof to
the classification of  all graphs $\graph$ such that
$$
\mbox{$\graph$ is geometric and $\ls|\graph|\geq49$};
$$
then, by
\eqref{eq-fano-must-be-hyperbolic},
the (plain) graph $\graph$ is a
$\fiber$-graph for a
certain affine Dynkin diagram \smash{$\fiber \ge \tA_2$}.
As in the proof of \autoref{thm-main}, we infer
that
$\graph$ is neither quadrangular (\autoref{prop.quad}), nor pentagonal
(\autoref{prop.pentagonal}),
nor astral (\autoref{prop.astral}), nor locally elliptic
(see \eqref{eq-loc-elliptic}).

For a geometric triangular graph $\graph$ with at least $49$
vertices and  $\rank(\graph) < 20$,
we apply \autoref{prop.smooth}
to show that $\graph$ is one of the five graphs \eqref{eq-monstrous-five}. Otherwise,
by \autoref{add.list.smooth}, the graph $\graph$
is one of the
rank~$20$
graphs that appear in \autoref{tab.smooth}.

To complete the deformation classification,
let $\graph$ be one of the graphs in \autoref{tab.smooth}.
As part of our study
of the saturation lists, we observe that the only geometric finite index
extension of $\Fano(\graph)$ is the trivial one; hence, one has
\[
\NS(X_4)=\Fano(\graph)\quad\text{and}\quad
\OG_h(\NS(X_4))=\Aut\graph
\label{eq.group}
\]
for any smooth quartic $X_4\subset\Cp2$ with $\Fn X_4\cong\graph$,
and the latter group is found using the \texttt{digraph}
package in \GAP~\cite{GAP4}. According to \cite[Theorem 3.9]{DIS},
the equilinear deformation families of
such quartics
are in a bijection with the primitive
isometric embeddings
\[
S:=\Fano(\graph)\into\bL=2\bE_8\oplus3\bU
\label{eq.embedding}
\]
regarded up to polarized autoisometry of $S\ni h$ and autoisometry of~$\bL$
preserving a coherent orientation of maximal positive definite subspaces
of~$\bL\otimes\R$ (the so-called \emph{positive sign structure});
such a family is
real if and only if \eqref{eq.embedding} admits a polarized autoisometry
reversing the positive sign structure. Hence, to complete the proof,
we classify embeddings~\eqref{eq.embedding}
using
Nikulin's~\cite{Nikulin:forms} theory of discriminant forms and
either Gauss~\cite{Gauss:Disquisitiones} theory of binary quadratic forms (in the
definite case $\rank T=2$) or Miranda--Morrison~\cite{Miranda.Morrison:book}
theory (in the indefinite case $\rank T\ge3$).
\qed

\remark
The groups $\Sym X_4$ and $\Aut(X_4,h)$ in \autoref{tab.smooth} are computed as
the subgroups of~\eqref{eq.group} that, respectively, act identically on
$\discr S$ or extend to an appropriate autoisometry of~$\bL$:
the latter is required to preserve the positive sign
structure (if $\rank T=2$) or act on~$T$ by $\pm1$ (if $\rank T\ge3$).
\endremark

\subsection{Proof of \autoref{add.real}}\label{proof.real}
As stated in \cite[Addendum 1.4]{DIS}, the 
number of real lines does take all
values in the range $\{0,\ldots,48\}$. Next, we recall that, when counting
the number of real lines on \emph{smooth} quartics, it suffices to consider
only those quartics~$X_4$ whose \emph{all} lines are real (with respect to
a certain real structure $\Gs\:X_4\to X_4$, see \cite[Proposition 3.10]{DIS}
or \cite[Theorem 2.7]{degt:lines}), and the latter is the case if and only if
the generic transcendental lattice~$T$ has a sublattice isomorphic to $[2]$
or $\bU(2)$, see \cite[Lemma 3.8]{DIS}. Hence, the statement of the addendum
follows from \autoref{tab.smooth} which lists all configurations of more than
$48$ lines and their respective transcendental lattices.
\qed

{
\let\.\DOTaccent
\def\cprime{$'$}
\bibliographystyle{amsplain}
\bibliography{degt-h5v6}
}

\end{document}

%% file: tab_smooth.tex
\config \bX_{64}\0(20)
 \ps{(6,0)^{16} (4,6)^{48}}&&&
  &\aut{4608} &\symplectic(192,1493) \group1152(6)
    &\rc{(1,0)} &\mat{[8,4,8]}\cr
\config \bX_{60}\i(20)
 \ps{(6,2)^{10} (4,4)^{30} (3,7)^{20}}&&&
  &\aut{480} &\symplectic(60,5) \group120(2)
    &\rc{(1,0)} &\mat{[4,2,16]}\cr
\config \bX_{60}\ii(20)
 \ps{(4,5)^{60}}&&&
  &\aut{240} &\symplectic(60,5) \group120(2)
    &\rc{(0,1)} &\mat{[4,1,14]}\cr
\config \bX_{56}\0(20)
 \ps{(4,6)^{8} (4,4)^{32} (2,8)^{16}}&&&
  &\aut{128} &\symplectic(16,11) \group64(4)
    &\rc{(0,1)} &\mat{[8,0,8]}\cr
\config \bY_{56}\0(20)
 \ps{(4,4)^{32} (3,7)^{24}}&&&
  &\aut{64} &\symplectic(16,8) \group32(2)
    &\rc{(1,0)} &\*\mat{[2,0,32]}\cr
\config \bQ_{56}\0(20)
 \ps{(4,4)^{24} (3,7)^{32}}&&&
  &\aut{384} &\symplectic(48,50) \group96(2)
    &\rc{(1,0)} &\mat{[4,2,16]}\cr
\config \bX_{54}\0(20)
 \ps{(6,2)^{4} (4,6)^{6} (4,4)^{6} (4,2)^{24} (2,8)^{12}}&&&
  &\aut{384} &\symplectic(24,12) \group48(2)
    &\rc{(1,0)} &\mat{[4,0,24]}\cr
\config \bQ_{54}\0(20)
 \ps{(4,4)^{24} (4,3)^{24}}&&&
  &\aut{48} &\symplectic(8,5) \group8(1)
    &\rc{(1,0)} &\mat{[4,2,20]}\cr
\config \bX_{52}\i(20)
 \ps{(6,0)^{1} (4,4)^{12} (4,3)^{12} (4,2)^{3} (3,5)^{18}}&&&
  &\aut{24} &\symplectic(3,1) \group3(1)
    &\rc{(1,0)} &\mat{[8,4,12]}\cr
\config \bX_{52}\ii(20)
 \ps{(6,0)^{1} (4,4)^{9} (4,3)^{18} (3,5)^{18}}&&&
  &\aut{36} &\symplectic(6,1) \group6(1)
    &\rc{(1,0)} &\mat{[4,2,20]}\cr
\config \bX_{52}\iii(20)
 \ps{(4,6)^{10} (3,5)^{40}}&&&
  &\aut{320} &\symplectic(20,3) \group80(4)
    &\rc{(1,0)} &\mat{[10,0,10]}\cr
\config \bX_{52}\v(20)
 \ps{(5,3)^{8} (3,5)^{32} (2,8)^{12}}&&&
  &\aut{32} &\symplectic(4,1) \group8(2)
    &\rc{(1,0)} &\mat{[10,4,10]}\cr
\config \bY_{52}\i(20)
 \ps{(4,6)^{2} (4,4)^{16} (3,5)^{20} (2,8)^{14}}&&&
  &\aut{8} &\symplectic(4,1) \group8(2)
    &\rc{(1,0)} &\*\mat{[2,0,38]}\cr&&&&&
    &\rc{(0,1)} &\mat{[8,2,10]}\cr
\config \bY_{52}\ii(20)
 \ps{(4,5)^{8} (4,3)^{12} (3,6)^{16} (2,7)^{16}}&&&
  &\aut{8} &\symplectic(4,2) \group8(2)
    &\rc{(1,0)} &\*\mat{[2,1,40]}\cr&&&&&
    &\rc{(0,1)} &\mat{[4,1,20]}\cr&&&&&
    &\rc{(0,1)} &\mat{[8,1,10]}\cr
\config \bZ_{52}\0(19)
 \ps{(6,0)^{4} (4,4)^{12} (4,2)^{24} (2,8)^{12}}&&&
  &\aut{384} &\symplectic(12,3) \group24(2)
    &\rc{(1,0)} &\*\mmat{\bU(2)\+[24]}\cr
\config \(52)(19)
 \ps{(4,3)^{48}}&&&
  &\aut{384} &\symplectic(8,5) \group8(1)
    &\rc{(1,0)} &\mmat{[-8]\+[4,2,4]}\cr
\config \bQ_{52}\i(20)
 \ps{(4,4)^{16} (4,3)^{16} (4,2)^{16}}&&&
  &\aut{64} &\symplectic(8,5) \group8(1)
    &\rc{(1,0)} &\mat{[4,0,24]}\cr
\config \bQ_{52}\ii(20)
 \ps{(4,4)^{8} (4,3)^{32} (4,2)^{8}}&&&
  &\aut{64} &\symplectic(16,11) \group16(1)
    &\rc{(0,1)} &\mat{[8,4,12]}\cr
\config \bQ_{52}\iii(20)
 \ps{(5,0)^{4} (3,6)^{48}}&&&
  &\aut{96} &\symplectic(4,1) \group24(6)
    &\rc{(1,0)} &\mat{[10,5,10]}\cr
\config \bX_{51}\0(20)
 \ps{(6,2)^{1} (5,3)^{6} (4,3)^{3} (3,6)^{6} (3,4)^{8} (2,7)^{27}}&&&
  &\aut{12} &\symplectic(6,1) \group6(1)
    &\rc{(0,1)} &\mat{[4,1,22]}\cr&&&&&
    &\rc{(1,0)} &\mat{[6,3,16]}\cr
\config \bX_{50}\i(20)
 \ps{(6,1)^{1} (4,4)^{9} (4,3)^{9} (4,2)^{9} (3,4)^{18}}&&&
  &\aut{18} &\symplectic(3,1) \group3(1)
    &\rc{(1,0)} &\mat{[4,2,28]}\cr
\config \bX_{50}\ii(20)
 \ps{(6,1)^{1} (4,4)^{6} (4,3)^{15} (4,2)^{6} (3,4)^{18}}&&&
  &\aut{12} &\symplectic(3,1) \group3(1)
    &\rc{(2,0)} &\mat{[4,0,24]}\cr
\config \bX_{50}\iii(20)
 \ps{(5,3)^{4} (4,4)^{8} (3,5)^{16} (2,8)^{4} (2,6)^{18}}&&&
  &\aut{16} &\symplectic(2,1) \group4(2)
    &\rc{(0,1)} &\mat{[4,0,24]}\cr
\config \bZ_{50}\0(19)
 \ps{(4,4)^{10} (3,5)^{40}}&&&
  &\aut{160} &\symplectic(10,1) \group20(2)
    &\rc{(1,0)} &\mmat{\bU(5)\+[4]}\cr
\config \(50)(19)
 \ps{(4,3)^{24} (4,2)^{24}}&&&
  &\aut{96} &\symplectic(8,5) \group8(1)
    &\rc{(1,0)} &\*\mmat{\bU(2)\+[28]}\cr
\config \bZ_{49}\0(19)
 \ps{(6,0)^{1} (4,3)^{18} (4,2)^{9} (3,4)^{18}}&&&
  &\aut{36} &\symplectic(3,1) \group3(1)
    &\rc{(1,0)} &\*\mmat{\bU(2)\+[28]}\cr